\documentclass{article}
\usepackage[utf8]{inputenc}
\usepackage[margin=1in]{geometry}
\usepackage{amsmath}
\usepackage{amsfonts}
\usepackage{graphicx}
\usepackage{algorithm,algcompatible}
\usepackage{overpic}
\usepackage{subcaption}
\usepackage{pict2e}
\usepackage{amsthm}
\usepackage[dvipsnames]{xcolor}
\usepackage{tabularx}
\usepackage{cite}
\usepackage{pgfplots}
\usepackage{multirow}

\usepackage[pdftex,
            pdfauthor={Rudy Geelen, Stephen Wright, Karen Willcox},
            pdftitle={Operator inference for non-intrusive model reduction with quadratic manifolds},
            pdfsubject={Data-driven model reduction, Nonlinear manifolds, Operator inference, Proper orthogonal decomposition},
            pdfkeywords={Data-driven model reduction, Nonlinear manifolds, Operator inference, Proper orthogonal decomposition}]{hyperref}
\hypersetup{colorlinks=true, linkcolor=blue, filecolor=magenta, urlcolor=cyan}

\usepackage{pgfplots}
\pgfplotsset{compat=1.15}
\usepgflibrary{fpu}
\usetikzlibrary{calc, colorbrewer, arrows}
\usepgfplotslibrary{fillbetween}

\DeclareMathOperator*{\argmin}{arg\,min}

\newtheorem{remark}{Remark}

\title{Operator inference for non-intrusive model reduction with quadratic manifolds}

\author{Rudy Geelen\thanks{Oden Institute for Computational Engineering and Sciences, University of Texas at Austin, Austin, TX (rudy.geelen@austin.utexas.edu, kwillcox@oden.utexas.edu, \url{https://kiwi.oden.utexas.edu})} \and Stephen Wright\footnotemark[2]\thanks{
Computer Sciences Department, University of Wisconsin, Madison, WI (swright@cs.wisc.edu)} \and Karen Willcox\footnotemark[1] }

\date{}
\begin{document}

\maketitle

\begin{abstract}
This paper proposes a novel approach for learning a data-driven quadratic manifold from high-dimensional data, then employing this quadratic manifold to derive efficient physics-based reduced-order models. The key ingredient of the approach is a polynomial mapping between high-dimensional states and a low-dimensional embedding. This mapping consists of two parts: a representation in a linear subspace (computed in this work using the proper orthogonal decomposition) and a quadratic component. The approach can be viewed as a form of data-driven closure modeling, since the quadratic component introduces directions into the approximation that lie in the orthogonal complement of the linear subspace, but without introducing any additional degrees of freedom to the low-dimensional representation. Combining the quadratic manifold approximation with the operator inference method for projection-based model reduction leads to a scalable non-intrusive approach for learning reduced-order models of dynamical systems. Applying the new approach to transport-dominated systems of partial differential equations illustrates the gains in efficiency that can be achieved over approximation in a linear subspace.
\end{abstract}

\section{Introduction}

Dimensionality reduction plays an important role in compressing high-dimensional datasets and in deriving reduced-order models that provide computationally efficient approximations of complex system dynamics.
Dimensionality reduction via projection onto a low-dimensional linear subspace underlies a large class of methods, including principal component analysis (PCA) and proper orthogonal decomposition (POD).
While linear dimension reduction is effective in a broad range of applications, many problems do not admit a fast decaying Kolmogorov \textit{n}-width, which limits the reduction that can be achieved using a linear subspace approximation \cite{pinkus2012n}.
To address this limitation, nonlinear dimensionality reduction techniques employ nonlinear manifolds, introducing richer structure to the low-dimensional representation.
In this paper, we propose a new approach for learning a quadratic manifold from high-dimensional data via the solution of a linear regression problem. We show how the quadratic manifold approximation can be combined with a non-intrusive model reduction approach to learn efficient physics-based reduced-order models of dynamical systems.
While the quadratic manifold approximation method proposed here has general applicability, our particular focus is on the reduced-order modeling of systems governed by parametrized partial differential equations (PDEs).

Our work is inspired by the quadratic manifold approaches developed for projection-based reduced-order models in structural dynamics problems featuring geometric nonlinearities \cite{JAIN201780, RUTZMOSER2017196,tatsis2018state}.
That work builds reduced-order models using quadratic manifolds comprised of vibration modes and modal derivatives, found through the solution of a generalized eigenvalue problem.
The manifold is tangent to a subspace spanned by the most relevant vibration modes, and its curvature is provided by modal derivatives obtained from sensitivity analysis of the eigenvalue problem.
Here, we seek a more general nonlinear manifold representation that can be employed to approximate the high-dimensional states that arise in parametrized PDEs.
The key ingredient in our work is a nonlinear mapping, postulated in polynomial form, which provides a basis for reduction.
The mapping can be constructed in non-intrusive fashion using linear regression and is driven by physics-based training data.
A similar concept was explored in \cite{https://doi.org/10.1002/nme.6831}, where a kernel principal component analysis is used to find a nonlinear manifold with lower dimensionality.
In that work, the approximation space is enriched with element-wise cross-products of the snapshots, thereby establishing globally curved manifolds.
Our approach differs from \cite{https://doi.org/10.1002/nme.6831} in that we explicitly derive a data-driven basis to represent the quadratic components of the manifold.

In this paper, we derive the approximation manifold from a representative set of snapshot data.
In the dynamical system setting, each snapshot represents the solution of the underlying governing equations at different points in time and/or different parameters. Snapshot-based approaches such as PCA and POD use the singular value decomposition of the snapshot matrix to compute a linear subspace defined by a set of basis vectors. The corresponding singular values quantify the dimension of the low-dimensional linear subspace (i.e., number of basis vectors) required to achieve a given level of accuracy in approximating the high-dimensional state snapshots.
Transport-dominated problems form one class of applications for which approximation in a linear subspace is typically inadequate, due to slow decay of the singular values.
Past approaches to address this challenge include transforming the basis to improve its approximation power.
Such methods include ``freezing'' \cite{OHLBERGER2013901}, shifting of the POD basis \cite{doi:10.1137/17M1140571}, and manifold calibration/transformation techniques \cite{PhysRevE.89.022923, Cagniart2019, https://doi.org/10.1002/nme.5998, doi:10.1137/19M1271270, refId0, doi:10.1137/16M1059904, rim2020manifold, Unger}.
These approaches typically rely on substantial additional knowledge about the underlying problem, such as the particular advection phenomena that govern the basis shifting.

Several approaches have been pursued to address the \textit{n}-width limitation of linear subspaces from a more general perspective.
These nonlinear model reduction methods aim to \textit{break} the Kolmogorov barrier by seeking approximations on nonlinear solution-manifolds instead of in a linear subspace.
Such constructions are characterized by a more rapid decay of the error with respect to number of degrees of freedom \cite{peherstorfer69breaking}. Some nonlinear methods rely on online adaptive model reduction \cite{https://doi.org/10.1002/nme.4800, doi:10.1137/19M1257275} and leverage the local nature of problems to derive efficient reduced models.
Other nonlinear model reduction techniques rely on the use of multiple linear subspaces (also \textit{dictionary} approaches) for constructing local approximation subspaces, instead of a single global approximation \cite{https://doi.org/10.1002/nme.4371, doi:10.1137/130924408, daniel2020model, geelen2021localized}.
More recently, efforts have been made to describe solution-manifolds through artificial neural networks \cite{LEE2020108973, KIM2021110841, GAO2020132614, DBLP:journals/corr/abs-2111-12995, SALVADOR20211, fresca2021comprehensive, KADEETHUM2022104098, doi:10.1063/5.0074310}.
These methods compute the nonlinear manifold in an offline stage and project the dynamical system onto this manifold in an online stage.
While these methods are becoming popular for reduced-order modeling of PDE models, questions about scalability and interpretability persist.

Our contribution in this paper is a novel reduced-order modeling framework for dynamical systems based on nonlinear solution-manifolds.
The nonlinear mapping is driven by physics-based training data.
It uses least-squares
(possibly regularized) to infer a quadratic mapping operator from the error that remains after projection onto a low-dimensional linear subspace of the full-order model.
The nonlinear manifold can thus be constructed efficiently in a non-intrusive manner from off-the-shelf solvers.
The manifold approach is then integrated with the operator inference method from \cite{PEHERSTORFER2016196} to infer reduced-order model operators from time domain simulation-data.
The overall strategy is a flexible framework for nonlinear reduction of low-order polynomial systems.
The non-intrusivity property of the proposed method does not imply a black-box formulation; rather, we use outputs of the full-order model, without having access to the code that produced the simulation data, and explicit knowledge in the form of a high-fidelity problem definition \cite{ghattas_willcox_2021}. 

Our methodology can be interpreted alternatively
as a \textit{closure} model, in that the reduced-order model of the original dynamical system is augmented with an additional term that accounts for the modeling error that remains after modal truncation.
The closure problem has a rich history in computational fluid dynamics and is typically
motivated and justified by invoking physical insights or mathematical arguments \cite{meneveau2006large}.
Closure modeling has been employed in the reduced-order modeling of complex PDE systems, where it emulates the effect of the discarded modes on the reduced model dynamics \cite{WANG201210, doi:10.1137/18M1177263, gouasmi2017priori, doi:10.1137/19M1292448}.
This property is particularly important because an indirect consequence of truncation is that the truncated models ignore the nonlinear interactions between the discarded and retained modes \cite{doi:10.1063/5.0061577}, so the inferred reduced models can fail to preserve all the solution features of interest in the dynamical system.
The formulation presented here explicitly accounts for the discarded modes through the inclusion of a data-driven closure term.
Importantly, this term is computed from observational data in a physics-agnostic fashion without requiring additional phenomenological arguments.

The paper is structured as follows.
In Section~\ref{sec:nonlinear_manifolds} we describe the data-driven construction of nonlinear solution-manifolds for dimensionality reduction, including an approach to reduce the number of basis functions used to define the nonlinear components of the manifold.
Section~\ref{sec:model_reduction} discusses how these nonlinear constructions can be leveraged to learn reduced-order models for high-dimensional dynamical systems.
In Section~\ref{sec:demonstration} we present numerical experiments for two dynamical systems arising from the discretization of linear PDEs.
Conclusions and proposed future developments are described in Section~\ref{sec:conclusions}.

\section{Data-driven learning of nonlinear manifolds for dimensionality reduction}
\label{sec:nonlinear_manifolds}

Dimensionality reduction seeks a low-dimensional representation of a high-dimensional dataset. In this section, we first review the widely used approach of representation in a low-dimensional linear subspace and then describe a nonlinear dimensionality reduction technique that uses a data-driven quadratic manifold.
To approximate the high-dimensional state $\mathbf{s}(t)\in \mathbb{R}^n$, we seek a mapping $\boldsymbol{\Gamma}: \mathbb{R}^r \mapsto \mathbb{R}^n$ with reduced dimension $r \ll n$ such that
\begin{equation}
    \mathbf{s}(t) \approx \boldsymbol{\Gamma}(\widehat{\mathbf{s}}(t)),
    \label{eq:nonlinear_mapping}
\end{equation}
where the vector $\widehat{\mathbf{s}}(t) \in \mathbb{R}^r$ denotes the reduced state coordinates of dimension $r$, and $t$ denotes some parameter on which the state depends.
The mapping $\boldsymbol{\Gamma}$ constitutes a transformation function that remains to be defined. Sections~\ref{subsec:linear} and \ref{subsec:nonlinear} will discuss explicit linear and nonlinear formulations for $\boldsymbol{\Gamma}$, respectively. An implementation of the quadratic manifold constructs and the numerical experiments carried out in this section are publicly available at \url{https://github.com/geelenr/quad_manifold}.


\subsection{The linear dimensionality reduction framework}
\label{subsec:linear}

Linear dimensionality reduction can be achieved via the singular value decomposition (SVD) using principal component analysis (PCA) \cite{pearson1901liii}) or (as it is known in the reduced-order modeling literature) the proper orthogonal decomposition (POD) \cite{lumley1967structure, sirovich1987turbulence, berkooz1993proper}.
We define the data matrix $\mathbf{S} \in \mathbb{R}^{n \times k}$, whose $j$th column is the state $\mathbf{s}(t_j):=\mathbf{s}_j$ (referred to as the $j$th snapshot), with a total of $k$ snapshots, where $k<n$.
Scaling of the data is important in obtaining an adequate basis.
To this end we introduce a reference matrix $\mathbf{S}_\text{ref}$, each of whose (identical) columns is reference state $\mathbf{s}_\text{ref}$ that shifts the training data. The reference state is chosen in a problem-specific manner.
The singular values of $\mathbf{S}-\mathbf{S}_\text{ref}$ are denoted by $\sigma_1 \geq \sigma_2 \geq \dots \geq \sigma_k \geq 0$.
The POD basis vectors are the left singular vectors of $\mathbf{S}-\mathbf{S}_\text{ref}$ corresponding to its $r$ largest singular values.
The original space is then reduced to the space spanned by an $r$-dimensional orthonormal set of basis vectors through the mapping
\begin{equation}
    \boldsymbol{\Gamma}(\widehat{\mathbf{s}}(t)) := \mathbf{s}_\text{ref} + \mathbf{V} \widehat{\mathbf{s}}(t),
    \label{eq:lin_subspace_approx}
\end{equation}
where $\mathbf{V}  = [\mathbf{v}_1 \,| \, \dotsc \,|\, \mathbf{v}_r] \in \mathbb{R}^{n \times r}$ is the POD basis matrix spanning a $r$-dimensional linear subspace to which the states are restricted, typically with $r \ll n$, and $\mathbf{v}_j \in \mathbb{R}^n$, $j=1,\ldots,r$, denotes the $j$th POD basis vector.

\begin{figure}[btp]
\centering \small
\begin{overpic}[width=0.6\linewidth]{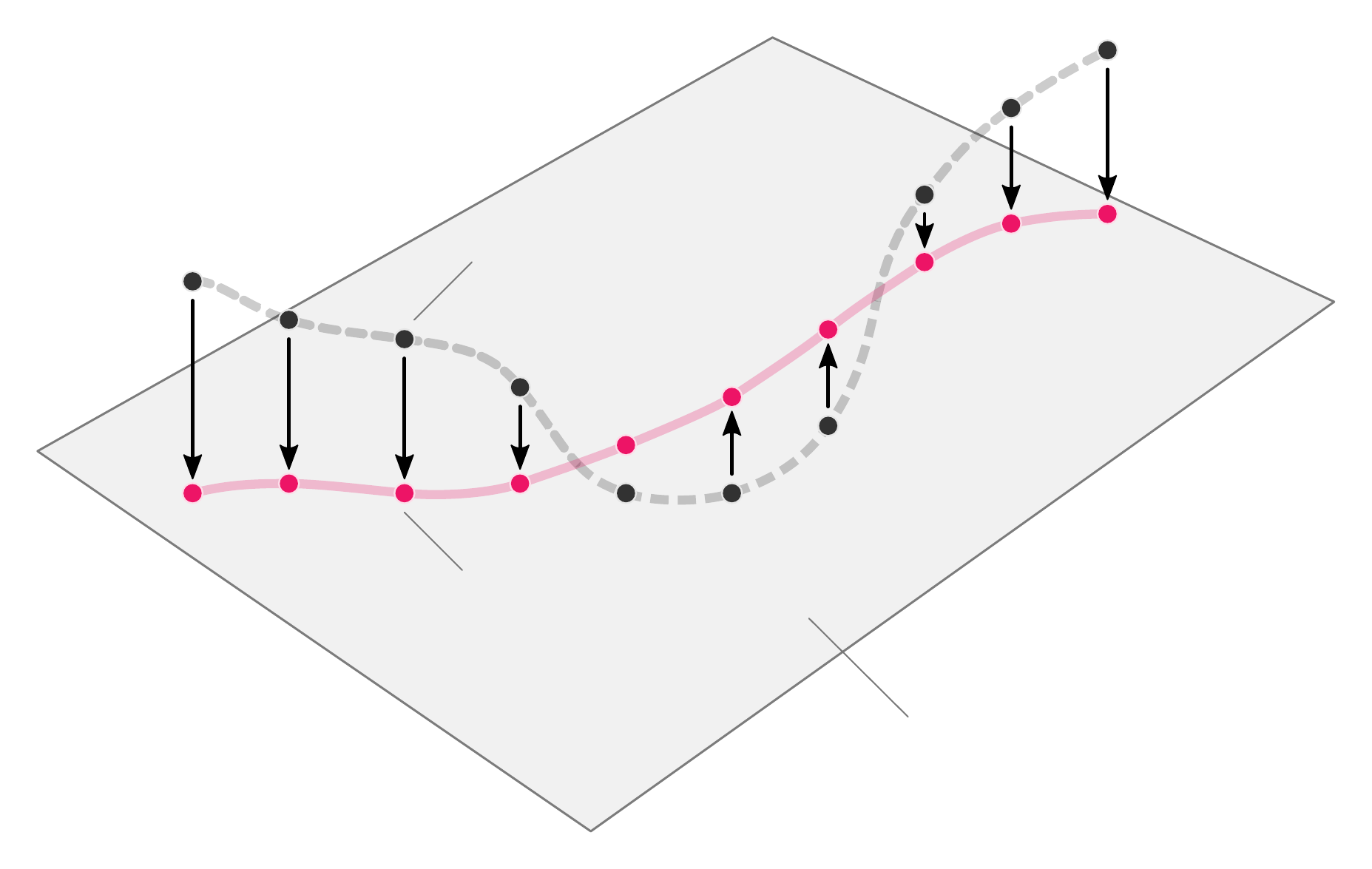}
\put(68,9){$\mathcal{V} = \text{range}(\mathbf{V})$}
\put(36,44){$\mathbf{s}_j-\mathbf{s}_\text{ref}$}
\put(27,17.5){\textcolor{WildStrawberry}{$\mathbf{V}\mathbf{V}^\top (\mathbf{s}_j-\mathbf{s}_\text{ref})$}}
\put(25,33){$\boldsymbol{\varepsilon}_j$}
\end{overpic}
\caption{Projection error using linear techniques for dimensionality reduction.}
\label{fig:projection_error}
\end{figure}
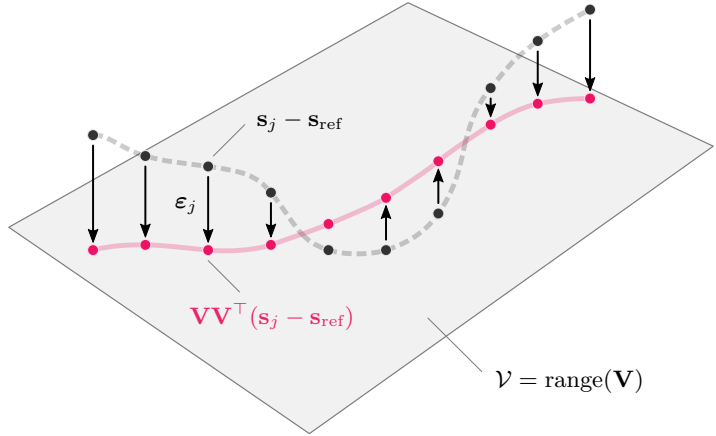

The POD basis minimizes the least-squares error of snapshot reconstruction. The familiar singular value decomposition result yields the projection error as being the sum of the squares of the singular values corresponding to those left singular vectors \textit{not} included in the POD basis:
\begin{equation}
\big\| \boldsymbol{\mathcal{E}} \big\|_F^2 = \sum_{j=r+1}^k \sigma_{j}^2,
\label{eq:snapshot_reconstruction_error}
\end{equation}
where
\begin{equation}
    \boldsymbol{\mathcal{E}} := (\mathbf{I} - \mathbf{V}\mathbf{V}^\top) (\mathbf{S}-\mathbf{S}_\text{ref}) := \mathbf{V}_\perp \mathbf{V}_\perp^\top (\mathbf{S}-\mathbf{S}_\text{ref}) \in \mathbb{R}^{n \times k}
    \label{eq:error}
\end{equation}
is the part of $\mathbf{S}-\mathbf{S}_\text{ref}$ not captured by the linear basis matrix $\mathbf{V}$.
In \eqref{eq:error}, $\mathbf{I}$ is the identity matrix of size $n$, and $\mathbf{V}_\perp \in \mathbb{R}^{n\times(n-r)}$ is the orthogonal complement of $\mathbf{V}$ in $\mathbb{R}^n$. The $j$th column of $\boldsymbol{\mathcal{E}}$, given by $\boldsymbol{\varepsilon}_j$, denotes the difference between the shifted snapshot $\mathbf{s}_j-\mathbf{s}_\text{ref}$ and its orthogonal projection onto a linear subspace $\mathcal{V}$, which is spanned by the columns of $\mathbf{V}$, as shown in Figure~\ref{fig:projection_error}.

In many datasets, driving the projection error \eqref{eq:snapshot_reconstruction_error} to an acceptably low value requires large $r$, rendering inefficient or impractical a reduction using linear techniques. In the following section we introduce a nonlinear dimensionality reduction approach that addresses this challenge by constructing a quadratic solution-manifold that incorporates directions in $\mathbf{V}_\perp$ in a data-driven fashion.

\subsection{Data-driven quadratic solution-manifolds}
\label{subsec:nonlinear}

In this section we derive a nonlinear manifold approach tailored for problems in which, for a desired reduced dimension $r$, the reconstruction error from \eqref{eq:snapshot_reconstruction_error} is non-negligible. We write
\begin{equation}
\mathbf{s}(t) = \mathbf{s}_\text{ref} + \mathbf{V} \widehat{\mathbf{s}}(t) + \boldsymbol{\varepsilon}(t) \in \mathbb{R}^n,
\label{eq:error_form}
\end{equation}
where $\boldsymbol{\varepsilon}(t)$ denotes the error associated with the projection of $\mathbf{s}(t)-\mathbf{s}_\text{ref}$ onto the linear subspace $\mathcal{V}$.
While many different manifold constructs are conceivable, we focus on a quadratic mapping between the high-dimensional data samples and their lower-dimensional representations.
Explicit nonlinear mappings with a polynomial structure were originally proposed in manifold learning \cite{6220279}, but similar formulations have emerged recently in model reduction \cite{JAIN201780, RUTZMOSER2017196, https://doi.org/10.1002/nme.6831, doi:10.1098/rspa.2021.0830, https://doi.org/10.48550/arxiv.2204.02462}.
Our attention will be restricted to the case of quadratic Kronecker products:
\begin{equation}
\boldsymbol{\varepsilon}(t)\approx \overline{\mathbf{V}} (\widehat{\mathbf{s}}(t) \otimes \widehat{\mathbf{s}}(t)),
\end{equation}
where $\overline{\mathbf{V}} \in \mathbb{R}^{n \times r^2}$ is a matricized quadratic mapping operator and $\otimes$ denotes the Kronecker product.\footnote{The operator $\otimes$ denotes the Kronecker product which for a column vector $\mathbf{s} = [s_1,s_2,...,s_n]^{\top}$ is given by $\mathbf{s} \otimes \mathbf{s} = \left[ s_1^2, s_1s_2, \dots, s_1s_n, s_2s_1, s_2^2, \dots, s_2s_n, \dots, s_n^2 \right]^\top \in \mathbb{R}^{n^2}.$}
The nonlinear transformation function is then
\begin{equation}
    \boldsymbol{\Gamma}(\widehat{\mathbf{s}}(t)) := \mathbf{s}_\text{ref} + \mathbf{V} \widehat{\mathbf{s}}(t) + \overline{\mathbf{V}} ( \widehat{\mathbf{s}}(t) \otimes \widehat{\mathbf{s}}(t) ).
    \label{eq:quadratic_approx}
\end{equation}
This type of nonlinearity yields data compression that improves on linear dimensionality reduction methods.
The approach readily generalizes to formulations with higher-order polynomial dependence.

The operator $\overline{\mathbf{V}}$ is inferred from the optimization problem
\begin{equation}
    \overline{\mathbf{V}} = \argmin_{\overline{\mathbf{V}} \in \mathbb{R}^{n \times r^2}} \sum_{j=1}^k \left\| \mathbf{s}_j - \mathbf{s}_\text{ref} - \mathbf{V} \widehat{\mathbf{s}}_j - \overline{\mathbf{V}} ( \widehat{\mathbf{s}}_j  \otimes \widehat{\mathbf{s}}_j ) \right\|_2^2,
    \label{eq:optim_problem}
\end{equation}
which means that the mapping-inferred operator satisfies the equations
\begin{equation}
\mathbf{s}_j \approx \mathbf{s}_\text{ref} + \mathbf{V} \widehat{\mathbf{s}}_j + \overline{\mathbf{V}} ( \widehat{\mathbf{s}}_j  \otimes \widehat{\mathbf{s}}_j ), \quad j=1,\dots,k,
\label{eq:approximate}
\end{equation}
at the optimum of the objective of \eqref{eq:optim_problem}.
In this representation, $\mathbf{V}$ is defined by the POD approach described  in the previous subsection and $\widehat{\mathbf{s}}_j$ is the reduced-order representation of the $j$th snapshot $\mathbf{s}_j$, defined by $\widehat{\mathbf{s}}_j = \mathbf{V}^\top (\mathbf{s}_j - \mathbf{s}_\text{ref})$.
By transposing the terms in the norm in the objective of \eqref{eq:optim_problem} and inserting the definition of the linear projection error, \eqref{eq:error}, the sum over the states can be written as a Frobenius norm:
\begin{equation}
    \overline{\mathbf{V}} = \argmin_{\overline{\mathbf{V}} \in \mathbb{R}^{n \times r^2}} \dfrac{1}{2}  \left\|  \mathbf{W}^\top \overline{\mathbf{V}} ^\top -  \boldsymbol{\mathcal{E}}^\top \right\|_F^2,
    \label{eq:least-squares}
\end{equation}
where $\boldsymbol{\mathcal{E}}$ is defined in \eqref{eq:error} and
\begin{equation}
    \mathbf{W} :=
    \begin{pmatrix}
    | & | & & | \\
    \widehat{\mathbf{s}}_1 \otimes \widehat{\mathbf{s}}_1 & \widehat{\mathbf{s}}_2 \otimes \widehat{\mathbf{s}}_2  & \dots & \widehat{\mathbf{s}}_k \otimes \widehat{\mathbf{s}}_k  \\
    | & | & & |
    \end{pmatrix} \in \mathbb{R}^{r^2 \times k}.
    \label{eq:quadratic_data}
\end{equation}
In the independent and simultaneous work \cite{https://doi.org/10.48550/arxiv.2204.02462}, a formulation that is mathematically identical to our Eq.~\eqref{eq:least-squares} is presented
and employed in a nonlinear model reduction approach. In the remainder of this section our work deviates from \cite{https://doi.org/10.48550/arxiv.2204.02462} in our analysis of the properties of $\overline{\mathbf{V}}$ and our numerical treatment of the optimization problem. Our deployment of the quadratic manifold approximation in reduced-order modeling also differs, as we discuss further in Section~\ref{sec:model_reduction}.

The optimization problem \eqref{eq:least-squares} finds the mapping operator $\overline{\mathbf{V}}$ that minimizes the error remaining after we use the columns of $\mathbf{W}$ to  fit the residual $\boldsymbol{\mathcal{E}}$ from the linear-manifold approximation, in a least-squares sense.
In essence, we \textit{learn} a more expressive mapping from a set of high-dimensional system-states to a low-dimensional subspace using the projection error $\boldsymbol{\mathcal{E}}$ as the driving mechanism. In matrix-vector format the terms in the norm in the objective of \eqref{eq:least-squares} can be written as
\begin{equation}
\mathbf{W}^\top \overline{\mathbf{V}} ^\top -  \boldsymbol{\mathcal{E}}^\top :=
\underbrace{\begin{pmatrix}
(\widehat{\mathbf{s}}_1 \otimes \widehat{\mathbf{s}}_1)^\top \\
(\widehat{\mathbf{s}}_2 \otimes \widehat{\mathbf{s}}_2)^\top \\
\vdots \\
(\widehat{\mathbf{s}}_k \otimes \widehat{\mathbf{s}}_k)^\top
\end{pmatrix}}_{k \times r^2}
\underbrace{\begin{pmatrix}
\overline{\mathbf{v}}_{1}^\top \\
\overline{\mathbf{v}}_{2}^\top \\
\vdots \\
\overline{\mathbf{v}}_{r^2}^\top
\end{pmatrix}}_{r^2 \times n}
- \underbrace{\begin{pmatrix}
\boldsymbol{\varepsilon}_1^\top \\
\boldsymbol{\varepsilon}_2^\top \\
\vdots \\
\boldsymbol{\varepsilon}_k^\top
\end{pmatrix}}_{k \times n},
\label{eq:ls_objective}
\end{equation}
from which it can be inferred that \eqref{eq:least-squares} is an overdetermined linear least-squares problem if $k > r^2$. By eliminating redundancy in the Kronecker products $(\widehat{\mathbf{s}}_j \otimes \widehat{\mathbf{s}}_j)$, we can reduce the number of columns in $\overline{\mathbf{V}}$ from $r^2$ to $r(r+1)/2$, and the condition for overdeterminedness becomes $k > r(r+1)/2$. Note that we can solve a vector least-squares problem for each row of $\overline{\mathbf{V}}$ in turn: a total of $n$ such problems, each with $r(r+1)/2$ unique reduced operator coefficients to be inferred.
Strictly speaking this means that we are working with a version of $\overline{\mathbf{V}}$ now that has dimension $n \times r(r+1)/2$, not $n \times r^2$ as in \eqref{eq:least-squares}--\eqref{eq:ls_objective}.
We can solve \eqref{eq:least-squares} explicitly via normal equations to obtain
\begin{equation}
\overline{\mathbf{V}}^\top = ( \mathbf{W} \mathbf{W}^\top)^{-1} \mathbf{W} \boldsymbol{\mathcal{E}}^\top
\quad \implies \quad
\overline{\mathbf{V}} = \boldsymbol{\mathcal{E}}  \mathbf{W}^\top ( \mathbf{W} \mathbf{W}^\top)^{-1} \in \mathbb{R}^{n \times r(r+1)/2}.
\label{eq:vbar_normal_equations}
\end{equation}
By combining this formula with the definition of $\boldsymbol{\mathcal{E}}$ from \eqref{eq:error}, we obtain
\begin{equation}
    \overline{\mathbf{V}} =
    (\mathbf{I}-\mathbf{V}\mathbf{V}^\top) (\mathbf{S}-\mathbf{S}_\text{ref}) \mathbf{W}^\top ( \mathbf{W} \mathbf{W}^\top)^{-1} =
    \mathbf{V}_\perp \left[ \mathbf{V}_\perp^\top (\mathbf{S}-\mathbf{S}_\text{ref}) \mathbf{W}^\top ( \mathbf{W} \mathbf{W}^\top)^{-1} \right] \in \mathbb{R}^{n \times r(r+1)/2},
    \label{eq:vbar_normal_equations_solution}
\end{equation}
from which it follows immediately that each column of $\overline{\mathbf{V}}$ is in the column space of $\mathbf{V}_\perp$, so that the orthogonality condition $\mathbf{V}^\top \overline{\mathbf{V}} = \mathbf{0}$ holds.

Ordinary least-squares estimators can suffer from noise amplification when the amount of training data is relatively small. This results in coefficients with large amplitudes, which are prone to overfitting the noise in the training set. A popular approach to avoid this problem is to add regularization to the least-squares cost function, which penalizes the norm of the coefficient vector. This promotes solutions that yield a good fit to the data using linear coefficients that are not too large. A commonly adopted strategy is Frobenius regularization, which replaces problem \eqref{eq:least-squares} by
\begin{equation}
    \overline{\mathbf{V}} := \argmin_{\overline{\mathbf{V}} \in \mathbb{R}^{n \times r(r+1)/2}}  \left( \dfrac{1}{2}  \left\|  \mathbf{W}^\top \overline{\mathbf{V}} ^\top -  \boldsymbol{\mathcal{E}}^\top \right\|_F^2 + \dfrac{\gamma}{2} \left\| \overline{\mathbf{V}} \right\|_F^2 \right),
    \label{eq:least-squares_reg1}
\end{equation}
where $\gamma$ is a scalar regularization parameter. This problem remains separable into $n$ vector least-squares problems. The explicit solution becomes
\begin{equation}
\overline{\mathbf{V}} = \boldsymbol{\mathcal{E}}  \mathbf{W}^\top ( \mathbf{W} \mathbf{W}^\top + \gamma \mathbf{I})^{-1} \in \mathbb{R}^{n \times r(r+1)/2},
\label{eq:vbar_normal_equations_reg}
\end{equation}
where the same argument yields that $\mathbf{V}^\top  \overline{\mathbf{V}} = \mathbf{0} $.

\begin{remark}[Simultaneous optimization of $\mathbf{V}$ and $\overline{\mathbf{V}}$]
The approach described above learns the quadratic mapping operator $\overline{\mathbf{V}}$ from the misfit of the linear dimension reduction with basis $\mathbf{V}$.
In theory, one could also pose the manifold learning problem to determine  the operators $\mathbf{V}$ and $\overline{\mathbf{V}}$ and the reduced-order representations $\widehat{\mathbf{s}}_j$, $j=1,2,\dotsc,k$ \textit{simultaneously}.
While such an approach would yield improved approximation accuracy, the formulation of this simultaneous determination is a  difficult constrained nonlinear optimization problem. Taking our cue from \eqref{eq:optim_problem},
the simultaneous optimization problem would be
\begin{equation}
    \min_{\mathbf{V}, \overline{\mathbf{V}}, \widehat{\mathbf{S}}} \, \sum_{j=1}^k \left\| \mathbf{s}_j - \mathbf{s}_\text{ref} - \mathbf{V} \widehat{\mathbf{s}}_j - \overline{\mathbf{V}} ( \widehat{\mathbf{s}}_j  \otimes \widehat{\mathbf{s}}_j ) \right\|_2^2,
    \label{eq:optim_problem2}
\end{equation}
where \[
    \widehat{\mathbf{S}} :=
    \begin{pmatrix}
    \widehat{\mathbf{s}}_1  & \widehat{\mathbf{s}}_2  & \dots & \widehat{\mathbf{s}}_k
    \end{pmatrix} \in \mathbb{R}^{r \times k}.
\]
If we retain the property of $\mathbf{V}$ having orthonormal columns, require  $\overline{\mathbf{V}}$ to satisfy $\mathbf{V}^\top \overline{\mathbf{V}}=\mathbf{0}$, and separate the minimization w.r.t.\ $\widehat{\mathbf{S}}$ from minimization w.r.t.\ $\mathbf{V}$ and $\overline{\mathbf{V}}$, we obtain
\begin{equation}
    \min_{\substack{\mathbf{V}, \overline{\mathbf{V}} \\ \mathbf{V}^\top\mathbf{V}=\mathbf{I} \\ \mathbf{V}^\top \overline{\mathbf{V}}=\mathbf{0}}} \min_{\widehat{\mathbf{S}}} \, \sum_{j=1}^k \left\| \mathbf{s}_j - \mathbf{s}_\text{ref} - \mathbf{V} \widehat{\mathbf{s}}_j - \overline{\mathbf{V}} ( \widehat{\mathbf{s}}_j  \otimes \widehat{\mathbf{s}}_j ) \right\|_2^2 =
     \min_{\substack{\mathbf{V}, \overline{\mathbf{V}} \\ \mathbf{V}^\top\mathbf{V}=\mathbf{I} \\ \mathbf{V}^\top \overline{\mathbf{V}}=\mathbf{0}}} \, \sum_{j=1}^k \min_{\widehat{\mathbf{s}}_j} \left\| \mathbf{s}_j - \mathbf{s}_\text{ref} - \mathbf{V} \widehat{\mathbf{s}}_j - \overline{\mathbf{V}} ( \widehat{\mathbf{s}}_j  \otimes \widehat{\mathbf{s}}_j ) \right\|_2^2.
    \label{eq:optim_problem3}
\end{equation}
By using the fact that for any vector $\mathbf{a}$ we have by orthonormality of $\mathbf{V}$ that $\| \mathbf{a} \|_2^2 = \| \mathbf{V}^\top \mathbf{a} \|_2^2 + \| (\mathbf{I} - \mathbf{V} \mathbf{V}^\top) \mathbf{a} \|_2^2$, and using $\mathbf{V}^\top\overline{\mathbf{V}}=\mathbf{0}$ we can rewrite  \eqref{eq:optim_problem3} as
\begin{equation}
    \min_{\substack{\mathbf{V}, \overline{\mathbf{V}} \\ \mathbf{V}^\top\mathbf{V}=\mathbf{I} \\ \mathbf{V}^\top \overline{\mathbf{V}}=\mathbf{0}}}  \, \sum_{j=1}^k \min_{\widehat{\mathbf{s}}_j} \Big\{ 
    \underbrace{\left\| \mathbf{V}^{\top} (\mathbf{s}_j - \mathbf{s}_\text{ref}) - \widehat{\mathbf{s}}_j \right\|_2^2}_\text{linear fit terms} +
    \underbrace{\left\| (\mathbf{I} - \mathbf{V} \mathbf{V}^\top) (\mathbf{s}_j - \mathbf{s}_\text{ref}) -  \overline{\mathbf{V}} ( \widehat{\mathbf{s}}_j  \otimes \widehat{\mathbf{s}}_j ) \right\|_2^2}_\text{quadratic fit terms} \Big\}.
    \label{eq:optim_problem4}
\end{equation}
In our proposed method, we obtain an approximate minimizer of \eqref{eq:optim_problem4} in two steps.
First, we choose $\mathbf{V}$ and $\widehat{\mathbf{s}}_j$, $j=1,2,\dotsc,k$ as the POD solution satisfying \eqref{eq:snapshot_reconstruction_error} and \eqref{eq:error}. 
(Note that this choice minimizes the sum of the linear fit terms in \eqref{eq:optim_problem4}, in fact, it makes each of these terms zero.)
Second, we substitute these values into the quadratic fit terms in \eqref{eq:optim_problem4}, and minimize the sum of these terms with respect to $\overline{\mathbf{V}}$,  yielding the formula \eqref{eq:vbar_normal_equations}, for which the constraint $\mathbf{V}^\top \overline{\mathbf{V}}=\mathbf{0}$ is satisfied.
Each of these steps is well defined and tractable, and yields a reasonable approximate solution to the problem \eqref{eq:optim_problem2}.
It remains an area of future work to assess the tradeoff between the potential improved approximation quality of joint optimization of $({\mathbf{V}}$, $\overline{\mathbf{V}}$, $\widehat{\mathbf{S}})$ and the increased computational complexity of solving this problem.
\end{remark}

\subsection{Data-driven learning of quasi-quadratic manifolds}
\label{subsec:quasi-quadratic}


The number of columns in $\overline{\mathbf{V}}$ scales as $\mathcal{O}(r^2)$.
Are all these basis functions necessary?
Could we span a nonlinear manifold adequately with a \textit{subset} of the basis functions that suffice to  define the nonlinear aspect of the manifold?
To this end, we explore the use of  \textit{column selection} to reduce the number of columns in  $\overline{\mathbf{V}}$.
We will use a well known technique for group-sparse regularization, treating each column as a ``group.''
The technique is to add a multiple of the sum-of-$\ell_2$ regularization function to the objective, with each term in the sum being the $\ell_2$ norm of a single column of $\overline{\mathbf{V}}$. We obtain
\begin{equation}
    \overline{\mathbf{V}} := \argmin_{\overline{\mathbf{V}} \in \mathbb{R}^{n \times r(r+1)/2}}  \left( \dfrac{1}{2}  \left\|  \mathbf{W}^\top \overline{\mathbf{V}} ^\top -  \boldsymbol{\mathcal{E}}^\top \right\|_F^2 + \gamma \sum_{j=1}^{r(r+1)/2} \left\| \overline{\mathbf{v}}_j \right\|_2 \right).
    \label{eq:least-squares_reg2}
\end{equation}
The sum-of-$\ell_2$ regularizer induces sparsity in a similar way to the $\ell_1$ norm of a vector, except that instead of producing element-wise sparsity, it induces sparsity {\em by groups}. By contrast, if we were to square each of the terms $\left\| \overline{\mathbf{v}}_j \right\|_2$, we would obtain the squared Frobenius norm of \eqref{eq:least-squares_reg1}, which does not sparsify $\overline{\mathbf{V}}$ at all.
An appropriate choice of the positive scalar parameter $\gamma$ forces several or many columns of $\overline{\mathbf{V}}$ to be zero at the minimizer of $\eqref{eq:least-squares_reg2}$.
Details of the algorithm for solving \eqref{eq:least-squares_reg2} are given in Appendix~\ref{app:sparse}.

Once we have selected the columns of $\overline{\mathbf{V}}$ (i.e.\ the rows of $\overline{\mathbf{V}}^\top$) that are nonzero at the solution of \eqref{eq:least-squares_reg2}, we perform a ``debiasing'' step, in which we solve a reduced version of \eqref{eq:least-squares} over just the columns of $\overline{\mathbf{V}}$ that were selected to be nonzero in \eqref{eq:least-squares_reg2}.
Since each column of $\overline{\mathbf{V}}$ corresponds to a row of $\mathbf{W}$, this reduced version of $\overline{\mathbf{V}}$ also involves a row submatrix of $\mathbf{W}$.
Essentially, we are replacing the quadratic mapping \eqref{eq:quadratic_approx} with the following \textit{quasi-quadratic} mapping:
\begin{equation}
    \boldsymbol{\Gamma}(\widehat{\mathbf{s}}(t)) := \mathbf{s}_\text{ref} + \mathbf{V} \widehat{\mathbf{s}}(t) + \overline{\mathbf{V}} ( \widehat{\mathbf{s}}(t)\, \widetilde{\otimes}\ \widehat{\mathbf{s}}(t) ),
    \label{eq:quadratic_approx2}
\end{equation}
where we have introduced the modified Kronecker product $\widetilde{\otimes}$ that performs the usual Kronecker multiplication, but then only takes a subset of size $q$ of its components, chosen, for example, by using the column selection procedure above.
This parametrization leads to a least-squares problem of the form \eqref{eq:least-squares}, but with a data matrix $\widetilde{\mathbf{W}} \in \mathbb{R}^{q \times k}$ modified accordingly, where $q<r(r+1)/2$ is a number of selected rows (the columns of $\widetilde{\mathbf{W}}^\top$)
\begin{equation}
    \widetilde{\mathbf{W}} :=
    \begin{pmatrix}
    | & | & & | \\
    \widehat{\mathbf{s}}_1 \, \widetilde{\otimes}\ \widehat{\mathbf{s}}_1 & \widehat{\mathbf{s}}_2 \, \widetilde{\otimes}\ \widehat{\mathbf{s}}_2  & \dots & \widehat{\mathbf{s}}_k \,\widetilde{\otimes}\ \widehat{\mathbf{s}}_k  \\
    | & | & & |
    \end{pmatrix} \in \mathbb{R}^{q \times k}.
\end{equation}
We continue to use $\overline{\mathbf{V}}$ to denote the matrix, now in $\mathbb{R}^{n \times q}$, that defines the nonlinear aspect of the solution-manifold. The terms in the norm in the objective of \eqref{eq:least-squares} can thus be rewritten as
\begin{equation}
\widetilde{\mathbf{W}}^\top \overline{\mathbf{V}} ^\top -  \boldsymbol{\mathcal{E}}^\top :=
\underbrace{\begin{pmatrix}
\widetilde{\mathbf{w}}_{1}^\top \\
\widetilde{\mathbf{w}}_{2}^\top \\
\vdots \\
\widetilde{\mathbf{w}}_{k}^\top
\end{pmatrix}}_{k \times q}
\underbrace{\begin{pmatrix}
\overline{\mathbf{v}}_{1}^\top \\
\overline{\mathbf{v}}_{2}^\top \\
\vdots \\
\overline{\mathbf{v}}_{q}^\top
\end{pmatrix}}_{q \times n} - \underbrace{\begin{pmatrix}
\boldsymbol{\varepsilon}_1^\top \\
\boldsymbol{\varepsilon}_2^\top \\
\vdots \\
\boldsymbol{\varepsilon}_k^\top
\end{pmatrix}}_{k \times n}.
\label{eq:ls_objective_efficient}
\end{equation}
The resulting problem is overdetermined if $k > q$ and can, once again, be decomposed into $n$ independent linear least-squares problems, each involving a vector of length $q$.
We note here that the earlier argument  based on normal equations suffices to show that the solution $\overline{\mathbf{V}}$ of \eqref{eq:ls_objective_efficient} will also satisfy the property $\mathbf{V}^\top \overline{\mathbf{V}}=\mathbf{0}$.

More information about the formulation \eqref{eq:least-squares_reg2}, including the proximal-gradient algorithm for solving this problem, the debiasing step, the sparsifying effect of the regularization term, and the choice of $\gamma$ can be found in \cite{wright2009sparse}. An earlier paper that used regularization like \eqref{eq:least-squares_reg2} is \cite{turlach2005simultaneous}.

\subsection{Constructing the quadratic manifold --- An illustrative 3D trajectory}
\label{subsec:model_problem}

The following toy example illustrates the construction of quadratic solution-manifolds via the proposed regression technique. Let $\mathbf{s}(t) \in \mathbb{R}^3$ be a given trajectory parametrized by the variable $t \in \mathcal{P}$, with $\mathcal{P} = [0,2\pi]$ a one-dimensional parameter space, as follows:
\begin{equation}
    \mathbf{s}(t) =
    \begin{pmatrix}
    s_1(t) \\
    s_2(t) \\
    s_3(t)
    \end{pmatrix} =
    \begin{pmatrix}
    \cos(t) \\
    \sin(t) \\
    \cos(2t)/2
    \end{pmatrix}.
    \label{eq:helix_trajectory}
\end{equation}
A dataset is built by uniformly sampling the trajectory $\mathbf{s}(t)$ at $k=100$ values of $t$. After computing the singular value decomposition of the shifted data matrix $\mathbf{S}-\mathbf{S}_\text{ref}$, where the columns of $\mathbf{S}_\text{ref}$ are given by the initial condition $\mathbf{s}(0)$, the left singular vectors are
\begin{equation}
\mathbf{v}_1 =
\begin{pmatrix}
   -0.9347\\
   \phantom{-}0\\
   -0.3554
\end{pmatrix}; \quad
\mathbf{v}_2 =
\begin{pmatrix}
   0  \\
   1 \\
   0
\end{pmatrix}; \quad
\mathbf{v}_3 =
\begin{pmatrix}
   \phantom{-}0.3554 \\
   \phantom{-}0 \\
   -0.9347
\end{pmatrix}.
\end{equation}
Representing the three-dimensional states in the POD coordinates $\widehat{\mathbf{s}}(t)$ can be achieved through a linear projection such as  $\widehat{\mathbf{s}}(t) = [\mathbf{v}_1, \mathbf{v}_2]^\top \mathbf{s}(t) \in \mathbb{R}^2$, where we have chosen a reduced dimension of $r=2$. The original trajectory $\mathbf{s}(t)$ is then approximated through $\mathbf{s}_\text{approx}(t) = \boldsymbol{\Gamma}(\widehat{\mathbf{s}}(t))$.

\begin{figure}[tb]
\centering
\begin{subfigure}[t]{.32\linewidth}
\begin{overpic}[width=\linewidth]{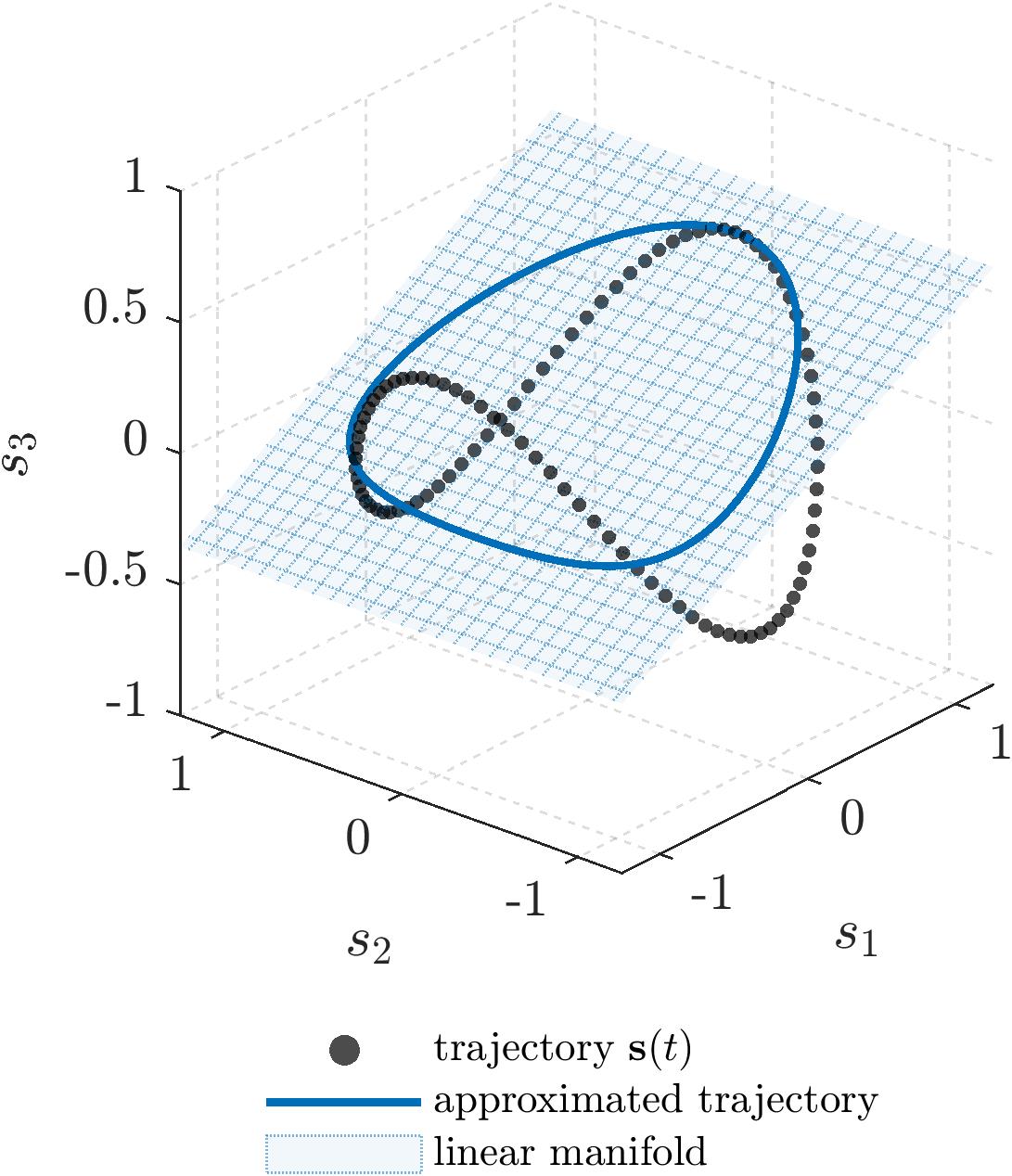}
\put(48,85){$\vector(-0.82,-1){5}$}
\put(48,85){$\vector(-1,0.37){7}$}
\put(47.5,81.25){\scriptsize $\mathbf{v}_1$}
\put(43.5,88){\scriptsize $\mathbf{v}_2$}
\end{overpic}
\caption{Linear subspace approach, \eqref{eq:lin_subspace_approx}.}
\label{fig:helix_linear}
\end{subfigure}
\begin{subfigure}[t]{.32\linewidth}
\begin{overpic}[width=\linewidth]{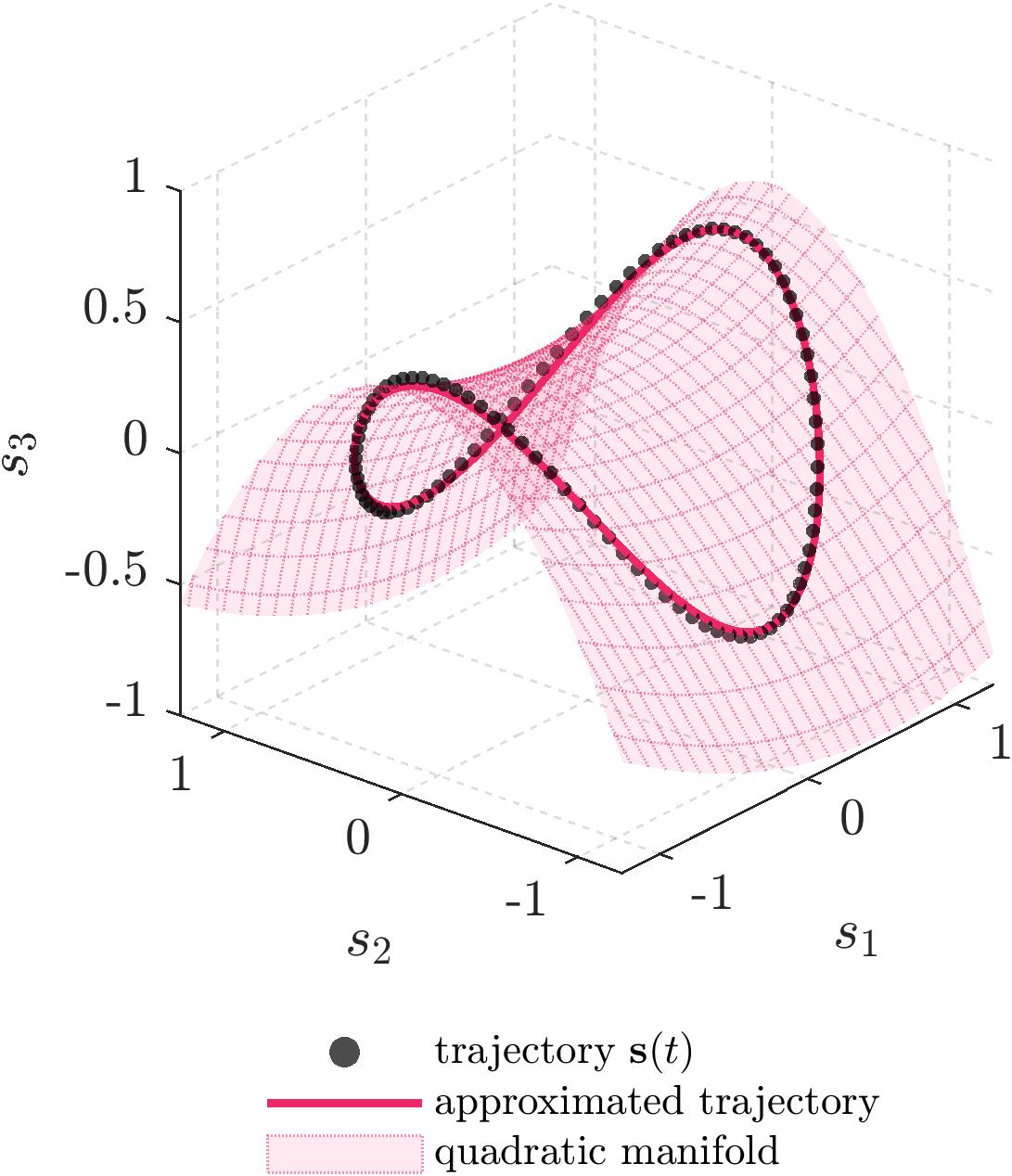}
\end{overpic}
\caption{Quadratic manifold approach, \eqref{eq:quadratic_approx};\eqref{eq:ls_objective} - no regularization.}
\label{fig:helix_quad}
\end{subfigure}
\begin{subfigure}[t]{.32\linewidth}
\begin{overpic}[width=\linewidth]{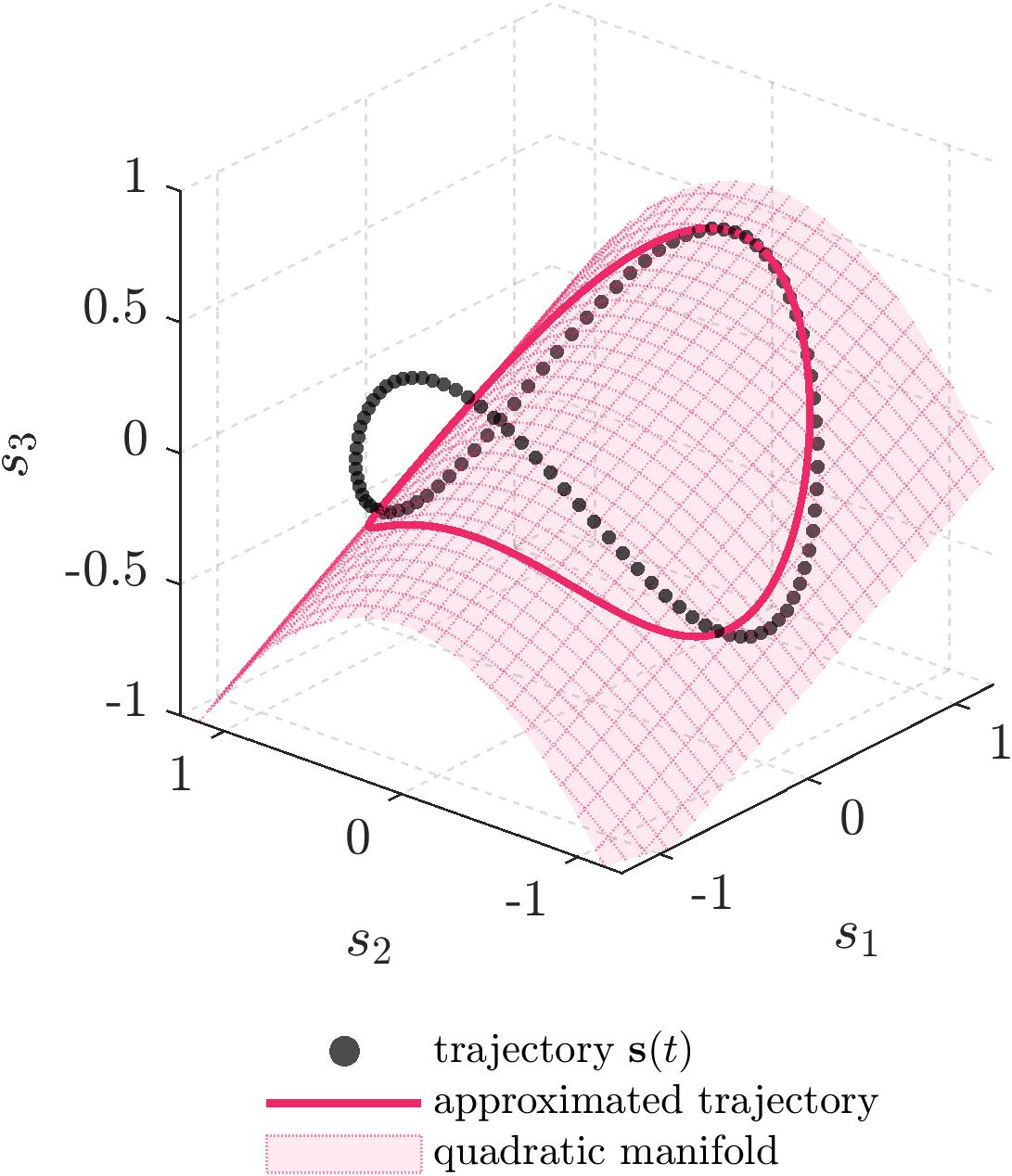}
\end{overpic}
\caption{Quasi-quadratic manifold approach, \eqref{eq:quadratic_approx};\eqref{eq:ls_objective_efficient} - no regularization.}
\label{fig:helix_data_efficient}
\end{subfigure}
\caption{Comparison of the linear and nonlinear manifold approaches. The black markers denote the sampled trajectory $\mathbf{s}(t)$, the solid lines correspond to their approximation $\mathbf{s}_\text{approx}(t)$, and the shaded surfaces denote the linear and nonlinear manifolds. The latter are visualized by means of projection of points in $\mathbb{R}^3$ onto the manifolds.}
\label{fig:helix}
\end{figure}

We now focus on the quadratic manifold formulation from Section \ref{subsec:nonlinear}, leading to least-squares problem \eqref{eq:least-squares} (i.e., with no regularization). Solving \eqref{eq:least-squares} with reduced space dimension $r=2$ we obtain the $r(r+1)/2=3$ basis vectors associated with the nonlinear part of the solution-manifold:
\begin{equation}
\overline{\mathbf{v}}_1 =
   \begin{pmatrix}
   -0.0681 \\
   \phantom{-}0 \\
   \phantom{-}0.1792
   \end{pmatrix}; \quad
\overline{\mathbf{v}}_2 =
   \begin{pmatrix}
   0 \\
   0 \\
   0
   \end{pmatrix}; \quad
\overline{\mathbf{v}}_3 =
   \begin{pmatrix}
   \phantom{-}0.3154 \\
   \phantom{-}0 \\
   -0.8296
   \end{pmatrix},
   \label{eq:vbar_cols}
\end{equation}
which correspond to the terms $\widehat{s}_1^2$, $\widehat{s}_1\widehat{s}_2$, and $\widehat{s}_2^2$, respectively. 
It can be seen that for this example, there exist only two basis vectors along which the linear subspace should be warped to optimally represent the training data in a least-squares sense. In this example, $\overline{\mathbf{v}}_2$ is the zero vector because including $\widehat{s}_1\widehat{s}_2$ in the quadratic manifold representation does not improve the approximation accuracy of the training data.
This is confirmed by the error metric from Table 1. The redundant basis vectors in such a construction will therefore be represented by zero-vectors.

We illustrate the quasi-quadratic nonlinear manifold formulation from Section~\ref{subsec:quasi-quadratic} for the scenario in which $q=1$ and we only employ the last column of $\mathbf{W}^\top$, corresponding to the term $\widehat{s}_2^2$. This column was identified by the column selection algorithm (Appendix \ref{app:sparse}) to be the most essential in spanning the quadratic solution-manifold. The basis matrix $\overline{\mathbf{V}}$ obtained by the procedure above consists of a single column, $\overline{\mathbf{v}}$ given by
\begin{equation}
\overline{\mathbf{v}} =
   \begin{pmatrix}
   \phantom{-}0.1638\\
   \phantom{-}0 \\
   -0.4308
   \end{pmatrix},
\end{equation}
which is \textit{not} among the set of basis vectors \eqref{eq:vbar_cols} computed without column selection, but has a similar direction to $\overline{\mathbf{v}}_3$.

Figure \ref{fig:helix} compares the reconstruction of the trajectory using the linear subspace approach (Section \ref{subsec:linear}), the quadratic manifold approach (Section \ref{subsec:nonlinear}), and the quasi-quadratic manifold approach (Section \ref{subsec:quasi-quadratic}). When the trajectory is restricted to the linear subspace spanned by the leading singular vectors $\mathbf{v}_1, \mathbf{v}_2$ (Figure~\ref{fig:helix_linear}), the approximate reconstruction of the full three-dimensional trajectory incurs a large projection error. The nonlinear manifold approaches induce a curving of the linear subspace without introducing additional degrees of freedom in the reduced-dimension representation. The full quadratic manifold provides an accurate approximation of the full trajectory (Figure~\ref{fig:helix_quad}). For the quasi-quadratic formulation we consider only one additional basis vector. This still leads to curving of the solution-manifold, but now along a single axis, as shown in Figure~\ref{fig:helix_data_efficient}. Table~\ref{tab:helix_error} lists the relative state prediction error for the different manifold approximations.
\begin{table}[htbp]
\centering \small
\caption{Comparison of the relative state prediction error for the linear subspace and quadratic manifold approaches. For the quadratic manifold approaches we indicate which of the columns of $\mathbf{W}^\top$ are used. We do not employ Frobenius regularization in this setting.}
\label{tab:helix_error}
\begin{tabular}{|l|l|l|}
\hline
\textbf{Reduction technique}  & \textbf{Selected columns} & \textbf{$\left. \| \mathbf{S} - \mathbf{S}_\text{approx} \|_F \middle/ \| \mathbf{S} \|_F \right.$} \\
\hline
\textbf{Linear manifold}  & /  & 0.4032 \\ \hline
\multirow{7}{*}{\textbf{Quadratic manifold}} & $\{1\}$ & 0.3971 \\
                                             & $\{2\}$ & 0.4032 \\
                                             & $\{3\}$ & 0.3042 \\
                                             & $\{1,2\}$ & 0.3971 \\
                                             & $\{2,3\}$ & 0.3042 \\
                                             & $\{1,3\}$ & 0.0258 \\
                                             & $\{1,2,3\}$ &  0.0258 \\ \hline
\end{tabular}
\end{table}

As expected, it can be seen that the linear subspace approach produces the largest error. The fully
quadratic manifold that takes into consideration all three columns of $\overline{\mathbf{V}}$ produces the lowest error.
As the number of columns in $\mathbf{W}^\top$ increases, the state prediction error drops and approaches the fully quadratic case where the data matrix is constructed using full Kronecker products. Because $\overline{\mathbf{v}}_2$ is a zero vector, it suffices to only consider the first and last columns. Indeed, the error in this case is the same as for the fully quadratic formulation. Furthermore, these results confirm that the third column (corresponding to $\widehat{s}_2^2$) is indeed the most essential in spanning the quadratic solution-manifold, as identified by the SpaRSA algorithm.

While strictly speaking regularization is not necessary if least-squares problem \eqref{eq:least-squares} is well-conditioned, we do advocate for its adoption: in practical problems we found that adequate regularization is a necessity in building stable and accurate data-driven reduced-order models, to be discussed next.


\section{Nonlinear manifolds for non-intrusive model reduction}
\label{sec:model_reduction}

We now show how the nonlinear manifold approximations introduced in Section~\ref{sec:nonlinear_manifolds} can be employed to derive reduced-order models for dynamical systems.
Section~\ref{subsec:proms_nonlinear} derives the form of a projection-based reduced model that results from the quadratic manifold approximation \eqref{eq:quadratic_approx}. Section~\ref{subsec:opinf} then proposes a quadratic manifold operator inference method in which the reduced models are learned from snapshot data in a non-intrusive fashion.

\subsection{Projection-based model reduction with nonlinear manifolds}
\label{subsec:proms_nonlinear}

Consider the high-dimensional dynamical system
\begin{equation}
\dfrac{\text{d}}{\text{d}t} \mathbf{s}(t) =  \mathbf{f}(t,\mathbf{s}(t)), \quad \mathbf{s}(0) = \mathbf{s}_0,
\label{eq:fom-nonlinear}
\end{equation}
where in this section $t\in (0,T]$ denotes time (with $T$ being the final time), the system state $\mathbf{s}(t) \in \mathbb{R}^n$ is time-dependent with specified initial condition $\mathbf{s}_0  \in \mathbb{R}^n$, and the right hand side is the nonlinear mapping $\mathbf{f}=\left[f_1, f_2, \dotsc, f_n\right]^\top : (0,T] \times  \mathbb{R}^n \rightarrow \mathbb{R}^n$. High-dimensional systems of the form \eqref{eq:fom-nonlinear} often result from discretization of partial differential equations. We will refer to \eqref{eq:fom-nonlinear} as the full-order model.

A reduced-order model of \eqref{eq:fom-nonlinear} is obtained by introducing the low-dimensional approximation \eqref{eq:nonlinear_mapping} and enforcing a Galerkin or Petrov-Galerkin orthogonality condition on the resulting residual.
Given a left-projection matrix $\boldsymbol{\Phi} \in \mathbb{R}^{n\times r}$, the reduced model is
\begin{equation}
\boldsymbol{\Phi}^\top \dfrac{\text{d}}{\text{d}t} \boldsymbol{\Gamma}(\widehat{\mathbf{s}}(t)) = \boldsymbol{\Phi}^\top \mathbf{f}(t,\boldsymbol{\Gamma}(\widehat{\mathbf{s}}(t))); \quad \widehat{\mathbf{s}}(0) = \widehat{\mathbf{s}}_0
\label{eq:fom_reduced-nonlinear}
\end{equation}
where $\widehat{\mathbf{s}}_0$ is the representation of the initial condition $\mathbf{s}_0$ in the reduced space coordinate system.
Consider linear dimensionality reduction where the basis $\mathbf{V}$ in \eqref{eq:lin_subspace_approx} is computed via  the POD as described in Section~\ref{subsec:linear}. Employing this linear approximation of the state together with a Galerkin projection (that is $\boldsymbol{\Phi}:=\mathbf{V}$) yields the reduced-order model
\begin{equation}
\dfrac{\text{d}}{\text{d}t} \widehat{\mathbf{s}}(t) = \mathbf{V}^\top \mathbf{f}(t,\mathbf{s}_\text{ref} + \mathbf{V} \widehat{\mathbf{s}}(t)); \quad \widehat{\mathbf{s}}(0) = \mathbf{V}^\top ( \mathbf{s}_0 - \mathbf{s}_\text{ref} ).
\label{eq:rom-nonlinear-linear}
\end{equation}
Consider instead quadratic dimensionality reduction \eqref{eq:quadratic_approx} with  $\mathbf{V}$ the POD basis and $\overline{\mathbf{V}}$ computed as in Section~\ref{subsec:nonlinear}. Again setting the left basis to be $\boldsymbol{\Phi}:=\mathbf{V}$,  we have the reduced-order model
\begin{equation}
\mathbf{V}^\top \left( \mathbf{V}\dfrac{\text{d}\widehat{\mathbf{s}}}{\text{d}t} + \overline{\mathbf{V}} \dfrac{\text{d}(\widehat{\mathbf{s}} \otimes \widehat{\mathbf{s}})}{\text{d}t} \right) = \dfrac{\text{d}}{\text{d}t} \widehat{\mathbf{s}}(t) = \mathbf{V}^\top \mathbf{f}(t,\mathbf{s}_\text{ref} + \mathbf{V} \widehat{\mathbf{s}}(t) + \overline{\mathbf{V}} ( \widehat{\mathbf{s}}(t) \otimes \widehat{\mathbf{s}}(t) )); \quad \widehat{\mathbf{s}}(0) = \mathbf{V}^\top ( \mathbf{s}_0 - \mathbf{s}_\text{ref} ),
\label{eq:rom-nonlinear-quad}
\end{equation}
where we have employed the orthogonality properties of the basis matrices (i.e., $\mathbf{V}^\top\mathbf{V}=\mathbf{I}$ and $\mathbf{V}^\top \overline{\mathbf{V}} = \mathbf{0}$).
Most existing model reduction methods employ the form \eqref{eq:rom-nonlinear-linear}; in \eqref{eq:rom-nonlinear-quad} we introduce the additional quadratic terms to the state approximation, however the dimension of the reduced-order model remains $r$ in both cases---that is, the quadratic manifold approximation does not introduce any additional degrees of freedom into $\widehat{\mathbf{s}}(t)$, although additional directions play a role in the state approximation through $\overline{\mathbf{V}}$.

Consider now a linear dynamical system, that is, $\mathbf{f}(t,\mathbf{s}(t))=\mathbf{A}\mathbf{s}(t)$, giving the full-order model
\begin{equation}
\dfrac{\text{d}}{\text{d}t} \mathbf{s}(t) =  \mathbf{A} \mathbf{s}(t); \quad \mathbf{s}(0) = \mathbf{s}_0,
\label{eq:fom}
\end{equation}
with $\mathbf{A} \in \mathbb{R}^{n \times n}$ the full-order linear operator. The Galerkin reduced-order model of \eqref{eq:fom} with linear approximation then takes the familiar form
\begin{equation}
\dfrac{\text{d}\widehat{\mathbf{s}}}{\text{d}t} = \widehat{\mathbf{c}} + \widehat{\mathbf{A}} \widehat{\mathbf{s}}; \quad \widehat{\mathbf{s}}(0) = \mathbf{V}^\top ( \mathbf{s}_0 - \mathbf{s}_\text{ref} ),
\label{eq:reduced_system_linear}
\end{equation}
where the reduced operators are $\widehat{\mathbf{c}} = \mathbf{V}^\top \mathbf{A} \mathbf{s}_\text{ref}$ and $\widehat{\mathbf{A}} = \mathbf{V}^\top \mathbf{A} \mathbf{V}$. Employing instead the quadratic approximation \eqref{eq:quadratic_approx}, the Galerkin reduced-order model of the linear system \eqref{eq:fom} takes the form
\begin{equation}
\mathbf{V}^\top \left( \mathbf{V}\dfrac{\text{d}\widehat{\mathbf{s}}}{\text{d}t} + \overline{\mathbf{V}} \dfrac{\text{d}(\widehat{\mathbf{s}} \otimes \widehat{\mathbf{s}})}{\text{d}t} \right) = 
\mathbf{V}^\top \mathbf{A} \mathbf{s}_\text{ref} + \mathbf{V}^\top \mathbf{A} \mathbf{V} \widehat{\mathbf{s}} + \mathbf{V}^\top \mathbf{A} \overline{\mathbf{V}}(\widehat{\mathbf{s}} \otimes \widehat{\mathbf{s}}); \quad \widehat{\mathbf{s}}(0) = \mathbf{V}^\top ( \mathbf{s}_0 - \mathbf{s}_\text{ref} ).
\label{eq:quadratic_full}
\end{equation}
Introducing the reduced matrix operators $\widehat{\mathbf{c}}, \widehat{\mathbf{A}}, \widehat{\mathbf{H}}$ as the constant, linear and quadratic components of the reduced-order model equation, and exploiting the orthogonality properties of the bases, \eqref{eq:quadratic_full} then simplifies to
\begin{equation}
\dfrac{\text{d}\widehat{\mathbf{s}}}{\text{d}t}  = \widehat{\mathbf{c}} + \widehat{\mathbf{A}} \widehat{\mathbf{s}} + \widehat{\mathbf{H}}(\widehat{\mathbf{s}} \otimes \widehat{\mathbf{s}}); \quad \widehat{\mathbf{s}}(0) = \mathbf{V}^\top ( \mathbf{s}_0 - \mathbf{s}_\text{ref} ),
\label{eq:reduced_system2}
\end{equation}
where
\begin{equation}
\widehat{\mathbf{c}} = \mathbf{V}^\top \mathbf{A} \mathbf{s}_\text{ref}, \quad \widehat{\mathbf{A}} = \mathbf{V}^\top \mathbf{A} \mathbf{V}, \quad
\widehat{\mathbf{H}} = \mathbf{V}^\top \mathbf{A} \overline{\mathbf{V}}.
\label{eq:reduced_order_operators}
\end{equation}

The quadratic manifold formulation transforms \textit{linear} problem \eqref{eq:fom} in the full state into a reduced system with \textit{quadratic} state dependence.
Correspondingly, if the full-order model had quadratic state dependency, approximation in the quadratic manifold would introduce a quartic state dependency in the reduced-order model.

The proposed approach may also be viewed as a form of data-driven closure modeling. In projection-based model reduction, closure modeling seeks to account for the effects of truncated modes. Our approach learns a data-driven closure term, $\overline{\mathbf{V}}(\widehat{\mathbf{s}} \otimes \widehat{\mathbf{s}})$ in \eqref{eq:quadratic_approx}, that introduces components in the range of $\mathbf{V}_{\perp}$ into the approximation of $\mathbf{s}$ (recall that $\mathbf{V}^\top \overline{\mathbf{V}} = \mathbf{0}$).
In our approach we do not explicitly account for the effect of the discarded POD modes on the retained ones (as in the scheme of \cite{doi:10.1137/19M1292448} for instance).
Instead, the unresolved dynamics are implicitly represented
via the quadratic term $\widehat{\mathbf{H}}(\widehat{\mathbf{s}} \otimes \widehat{\mathbf{s}})$ in the reduced state dynamical system \eqref{eq:reduced_system2}.

\begin{remark}[Quasi-quadratic formulation]
Using the quasi-quadratic formulation from Section~\ref{subsec:quasi-quadratic} for inferring the mapping operator $\overline{\mathbf{V}}$ leads to the modified reduced model
\begin{equation}
\dfrac{\mbox{\rm d}\widehat{\mathbf{s}}}{\mbox{\rm d}t}  = \widehat{\mathbf{c}} + \widehat{\mathbf{A}} \widehat{\mathbf{s}} + \widetilde{\mathbf{H}}(\widehat{\mathbf{s}}\: \widetilde{\otimes}\: \widehat{\mathbf{s}}); \quad \widehat{\mathbf{s}}(0) = \mathbf{V}^\top ( \mathbf{s}_0 - \mathbf{s}_\text{ref} ),
\label{eq:reduced_system3}
\end{equation}
where $\widetilde{\mathbf{H}} = \mathbf{V}^\top \mathbf{A} \overline{\mathbf{V}} \in \mathbb{R}^{r\times q}$ with $q<r(r+1)/2$.
This leads to a more data-efficient model to learn, but note that \eqref{eq:reduced_system_linear}, \eqref{eq:reduced_system2} and \eqref{eq:reduced_system3} all have reduced state dimension $r$. In the next subsection we describe construction of \eqref{eq:reduced_system2}, but the algorithm applies to constructing \eqref{eq:reduced_system3}, interchanging $\widetilde{\otimes}$ with $\otimes$ and $\widetilde{\mathbf{H}}$ with $\widehat{\mathbf{H}}$.
\end{remark}

\begin{remark}[Computational cost]
While the quadratic (or quasi-quadratic) approximation does not introduce additional degrees of freedom to the reduced state, it does incur additional computational cost. Consider a linear mapping of the form \eqref{eq:lin_subspace_approx} with reduced state dimension $r_\ell$, and consider a quasi-quadratic mapping of the form (\ref{eq:quadratic_approx2}) with reduced state dimension $r_q$ and $q$ columns in $\overline{\mathbf{V}}$. 
The reduced operators $\widehat{\mathbf{c}}$, $\widehat{\mathbf{A}}$ and $\widetilde{\mathbf{H}}$ in the quadratic manifold reduced model (\ref{eq:reduced_system3}) require storage of $r_q+r_q^2+qr_q$ elements. 
By contrast, the linear subspace reduced-order model (\ref{eq:reduced_system_linear}) requires storage of $r_\ell+r_\ell^2$ elements. 
Solving the linear reduced-order model (\ref{eq:reduced_system_linear}) requires evaluating at each timestep the matrix vector product $\widehat{\mathbf{A}}\widehat{\mathbf{s}}$, which entails $2r_\ell^2$ FLOPS. Solving the quadratic reduced-order model (\ref{eq:reduced_system3}) requires evaluating at each timestep the matrix vector product $\widehat{\mathbf{A}}\widehat{\mathbf{s}}$ ($2r_q^2$ FLOPS), evaluating $\widehat{\mathbf{s}}\: \widetilde{\otimes}\: \widehat{\mathbf{s}}$ ($q$ FLOPS), and evaluating  $\widetilde{\mathbf{H}}(\widehat{\mathbf{s}}\: \widetilde{\otimes}\: \widehat{\mathbf{s}})$ ($2qr_q$ FLOPS). For a given target accuracy in representing the snapshot data, we expect $r_q<r_\ell$, since the quadratic manifold representation achieves a given level of representation accuracy with a lower dimension. Whether the quadratic manifold reduced model is cheaper to solve than the linear subspace reduced model will be problem specific, depending on the compression from $r_\ell$ to $r_q$. It will also depend on the particular problem structure and whether, for example, it might be possible to factor $\widetilde{\mathbf{H}}$ in order to reduce online evaluation costs. Consideration of computational costs highlights the importance of the quasi-quadratic formulation: with the full quadratic approximation, we have $q=r_q(r_q+1)/2$ and the worst-case cost of evaluating $\widetilde{\mathbf{H}}(\widehat{\mathbf{s}}\: \widetilde{\otimes}\: \widehat{\mathbf{s}})$ scales with $r_q^3$. In this case, the full quadratic manifold reduced models are likely only to be competitive in cost for small $r_q$; as $r_q$ grows their cost will likely exceed that of a linear subspace reduced model that achieves the same level of accuracy with $r_\ell>r_q$.
\end{remark}

\begin{remark}[Initial condition satisfaction]
Satisfying the specified initial condition of the dynamical system can be important for achieving an accurate approximation in a nonlinear manifold \cite{LEE2020108973}. This requires the initial reduced coordinates $\widehat{\mathbf{s}}(0) := \widehat{\mathbf{s}}_0$ to satisfy $\boldsymbol{\Gamma}(\widehat{\mathbf{s}}_0) = \mathbf{s}_0$. In the case of quadratic solution-manifolds, as described in Section~\ref{sec:nonlinear_manifolds}, this implies
\begin{equation}
    \mathbf{s}_\text{ref} := \mathbf{s}_0 - \mathbf{V}\widehat{\mathbf{s}}_0 - \overline{\mathbf{V}}(\widehat{\mathbf{s}}_0 \otimes \widehat{\mathbf{s}}_0).
\end{equation}
It can be shown in a straightforward manner that this condition is satisfied exactly if $\mathbf{s}_\text{ref} = \mathbf{s}_0$. When dealing with trajectories computed from simulations with \textit{multiple} different initial conditions, one of these has to be chosen as the reference state, $\mathbf{s}_\text{ref}$, and exact satisfaction of initial conditions can no longer be guaranteed in all cases.
\end{remark}


\subsection{Learning quadratic manifold reduced models with operator inference}
\label{subsec:opinf}

Eq.~\eqref{eq:rom-nonlinear-quad} defines the quadratic manifold reduced model in the general nonlinear case and \eqref{eq:reduced_system2} defines the quadratic manifold reduced model in the case of a linear full-order model. While in some applications it may be possible to construct these reduced models intrusively (i.e., by explicitly computing the projected operators $\mathbf{V}^\top \mathbf{f}$, $\mathbf{V}^\top \mathbf{A}$, $\mathbf{V}^\top \mathbf{A} \mathbf{V}$ and $\mathbf{V}^\top \mathbf{A} \overline{\mathbf{V}}$), introduction of the quadratic manifold approximation complicates the implementation of an already time-consuming process. Eq.~\eqref{eq:reduced_system2} also illustrates that the quadratic manifold approximation changes the form of the reduced-order model, which further complicates implementation. We therefore employ the non-intrusive operator inference method of \cite{PEHERSTORFER2016196}, which infers the reduced model operators directly from time-domain simulation data.

A non-intrusive model reduction method is defined in \cite{ghattas_willcox_2021} as a method that computes the reduced model using outputs of the full-order model without having access to the full-order operators (or to their action on a vector). Non-intrusive methods are not necessarily black-box; they can exploit knowledge of the full-order problem definition and structure (in our case, the structure of the PDEs that govern the problem of interest). Non-intrusivity is important for cases where modifying the full-order model source code that produces the simulation data may be impractical or its internal access is restricted (e.g., legacy or commercial codes). In our case, non-intrusivity simplifies the numerical process to compute the quadratic manifold reduced-order model (\ref{eq:reduced_system2}).
Our use of operator inference with a quadratic manifold approximation differs from past contributions such as \cite{PEHERSTORFER2016196, QIAN2020132401, doi:10.1080/03036758.2020.1863237, BENNER2020113433}, which all employ operator inference with linear dimensionality reduction.

Given a data set comprising $k$ state snapshots $\mathbf{s}_1,\mathbf{s}_2,\dots,\mathbf{s}_k$ and the corresponding time derivative data, operator inference generates {\em reduced state} snapshots via projection with the basis matrix  $\mathbf{V}$, giving $\widehat{\mathbf{s}}_j=\mathbf{V}^\top (\mathbf{s}_j-\mathbf{s}_\text{ref}); j=1,\dots,k$.
For the linear full-model form \eqref{eq:fom}, operator inference then finds the quantities $\widehat{\mathbf{c}}$ and $\widehat{\mathbf{A}}$ that define the reduced model that best matches these projected snapshot data in a minimum residual sense, as follows:
\begin{equation}
    \left(\widehat{\mathbf{c}}, \widehat{\mathbf{A}}\right) = \argmin_{ \widehat{\mathbf{c}} \in \mathbb{R}^r, \widehat{\mathbf{A}} \in \mathbb{R}^{r \times r} } \sum_{j=1}^k
 \left\| \widehat{\mathbf{c}} + \widehat{\mathbf{A}} \widehat{\mathbf{s}}_j - \dfrac{\text{d}\widehat{\mathbf{s}}_j}{\text{d}t} \right\|_2^2  + \lambda \left( \| \widehat{\mathbf{c}} \|_2^2 + \| \widehat{\mathbf{A}} \|_F^2 \right),
 \label{eq:opinf_linear}
\end{equation}
where the second term is a Tikhonov regularization term and $\lambda>0$ is a scalar regularization parameter.
The time derivatives must be estimated, for example using a finite difference approximation.
As shown in \cite{PEHERSTORFER2016196}, this optimization problem decouples into $r$ independent linear least-squares problems, each involving a single row of $[\widehat{\mathbf{c}}\ \widehat{\mathbf{A}}]$.
The Tikhonov regularization term in \eqref{eq:opinf_linear} promotes stability and accuracy of solutions to the operator inference problem and inhibits overfitting of the system operators to the data \cite{doi:10.1080/03036758.2020.1863237} in the presence of noise, which can arise from error in the numerically estimated time derivatives, model misspecification, and unresolved system dynamics.

Inferring the quadratic manifold reduced model \eqref{eq:reduced_system2} is straightforward in the operator inference framework. We simply include the quadratic operator in the inference problem, to obtain
\begin{equation}
    \left(\widehat{\mathbf{c}}, \widehat{\mathbf{A}}, \widehat{\mathbf{H}} \right) = \argmin_{\substack{\widehat{\mathbf{c}} \in \mathbb{R}^r, \widehat{\mathbf{A}} \in \mathbb{R}^{r \times r}, \\ \widehat{\mathbf{H}} \in \mathbb{R}^{r \times r^2}}} \sum_{j=1}^k
 \left\|  \widehat{\mathbf{c}} + \widehat{\mathbf{A}} \widehat{\mathbf{s}}_j + \widehat{\mathbf{H}}(\widehat{\mathbf{s}}_j \otimes \widehat{\mathbf{s}}_j) - \dfrac{\text{d}\widehat{\mathbf{s}}_j}{\text{d}t} \right\|_2^2 + \lambda_1 \left( \| \widehat{\mathbf{c}} \|_2^2 + \| \widehat{\mathbf{A}} \|_F^2 \right)+ \lambda_2 \| \widehat{\mathbf{H}} \|_F^2\ ,
 \label{eq:opinf_quadratic}
 \end{equation}
where $\lambda_1$ and $\lambda_2$ are scalar regularization parameters (we typically choose to weight regularization of $\widehat{\mathbf{A}}$ and $\widehat{\mathbf{H}}$ differently as they are characterized by different scales). We also exploit the symmetry of $\widehat{\mathbf{H}}$ to reduce the number of inferred operator coefficients by not solving for its redundant terms. While the states remain $r$-dimensional vectors, the total number of coefficients of the unknown reduced operators to be inferred in \eqref{eq:opinf_quadratic} is greater than in \eqref{eq:opinf_linear}.

As noted above, if the full-order model had quadratic form, the quadratic manifold reduced-order model would have quartic structure. It is straightforward to formulate the inference of a quartic model in an analogous way to \eqref{eq:opinf_quadratic} (noting that the operator inference minimization remains a linear least-squares problem); however, the number of operator coefficients to be inferred grows rapidly. Reduced-order models with quadratic and quartic structure have $\mathcal{O}(r^2)$ and $\mathcal{O}(r^4)$ reduced operator coefficients to be inferred, respectively, for each of the $r$ independent least-squares problems. Even with column sub-selection, this may cause the operator inference problem to be under-determined and/or ill-conditioned.

Algorithm~\ref{alg:opinf_manifold} summarizes the complete approach for deriving a quadratic manifold reduced model of a linear dynamical system. We emphasize that the algorithm requires only a set of state snapshots and knowledge of the structure of the full-order model, not access to the full-order model code itself.

\begin{algorithm}
\caption{Quadratic Manifold Operator Inference.}
\label{alg:opinf_manifold}
\begin{algorithmic}[1]
\REQUIRE Snapshot matrix $\mathbf{S} \in \mathbb{R}^{n \times k}$, snapshot time derivatives $\dot{\mathbf{S}} \in \mathbb{R}^{n \times k}$, user-defined tolerance $\kappa$
\ENSURE Reduced model operators $\widehat{\mathbf{c}},\widehat{\mathbf{A}}, \widehat{\mathbf{H}}$, basis matrices $\mathbf{V}, \overline{\mathbf{V}}$
\item[]
\STATEx \COMMENT{\textit{Computing the linear POD basis}}
\STATE Center the snapshot data about a reference state
\STATE Compute the SVD of $\mathbf{S}-\mathbf{S}_\text{ref}$
\STATE $r \leftarrow $ Determine linear POD basis dimension (informed by the singular value decay and tolerance $\kappa$)
\STATE $\mathbf{V} \leftarrow $ The $r$ leading left singular vectors of $\mathbf{S}-\mathbf{S}_\text{ref}$
\item[]
\STATEx \COMMENT{\textit{Computing the quadratic manifold basis}}
\STATE Represent states in POD coordinates: $\widehat{\mathbf{S}} \leftarrow \mathbf{V}^\top (\mathbf{S}-\mathbf{S}_\text{ref})$
\STATE $\boldsymbol{\mathcal{E}}, \mathbf{W} \leftarrow $ Compute linear projection error \eqref{eq:error} and quadratic data matrix \eqref{eq:quadratic_data}
\IF {quasi-quadratic manifold formulation}
\STATE Column selection using the SpaRSA algorithm \eqref{alg:sparsa}: $\widetilde{\mathbf{W}} \leftarrow \mathbf{W}$
\ENDIF
\STATE $\overline{\mathbf{V}} \leftarrow $ Solve regularized least-squares problem using \eqref{eq:least-squares_reg1} or \eqref{eq:least-squares_reg2}
\item[]
\STATEx \COMMENT{\textit{Operator inference for learning low-dimensional dynamical systems}}
\STATE Represent time derivatives in POD coordinates: $\dot{\widehat{\mathbf{S}}} \leftarrow \mathbf{V}^\top \dot{\mathbf{S}}$
\STATE ${\lambda_1,\lambda_2} \leftarrow $ Set the regularization parameters
\STATE $\widehat{\mathbf{c}},\widehat{\mathbf{A}},\widehat{\mathbf{H}} \leftarrow $ Solve regularized operator inference problem \eqref{eq:opinf_quadratic}
\end{algorithmic}
\end{algorithm}

\section{Demonstration of the approach}
\label{sec:demonstration}

In this section we demonstrate the proposed quadratic manifold approach on two numerical examples: the advection equation and the wave equation. The advection equation describes the transport of a substance or quantity through advection, and is therefore of fundamental importance in physics and engineering sciences. The wave equation arises in fields such as acoustics, electromagnetism, and fluid dynamics, and is commonly used in studying various types of mechanical and electromagnetic waves. Both are transport-dominated problems for which the Kolmogorov $n$-width is well known to be problematic for achieving efficient model reduction in a static linear subspace \cite{doi:10.1137/19M1257275}. We calibrate the hyper-parameters $\lambda_1, \lambda_2, \gamma$ by choosing them to  minimize the relative error between the reduced-order model predictions and the available training data,
as in \cite{doi:10.1080/03036758.2020.1863237}.

\begin{figure}[tp]
\centering \small
\begin{subfigure}[c]{.49\linewidth}
\includegraphics[scale=1]{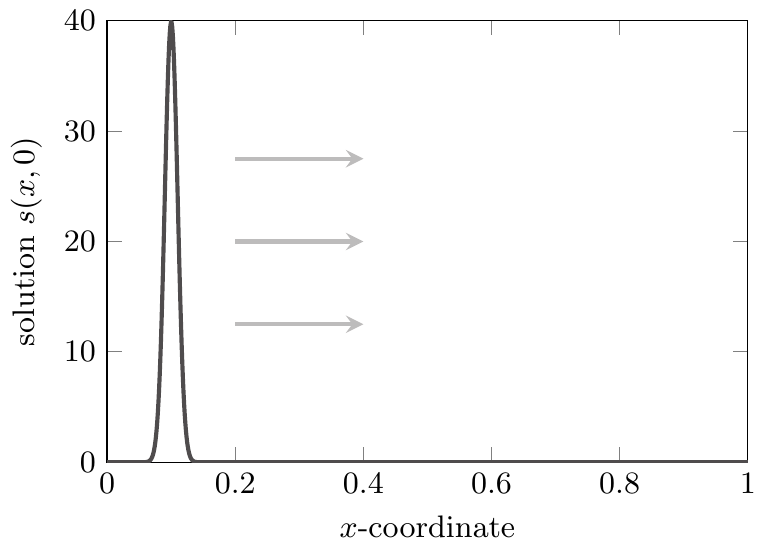}
\caption{Initial condition at $\mu=0.10$ and transport direction.}
\label{fig:1d_ic}
\end{subfigure}
\begin{subfigure}[c]{.49\linewidth}
\includegraphics[scale=1]{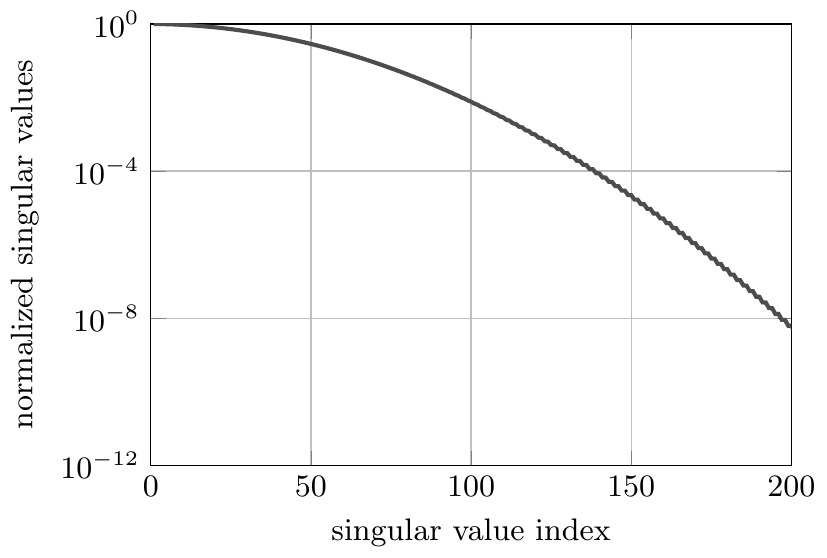}
\caption{Singular values.}
\label{fig:1d_svs}
\end{subfigure}
\caption{The one-dimensional advection equation with $c=10$. Plot (a) shows the initial condition and the direction of the transport. Plot (b) shows the decay of the singular values of the snapshot matrix spanning all the initial conditions.}
\label{fig:1d}
\end{figure}

\subsection{Linear transport equation}


\begin{figure}[tbp]
\centering \small
\begin{subfigure}[t]{.6\linewidth}
\includegraphics[scale=1]{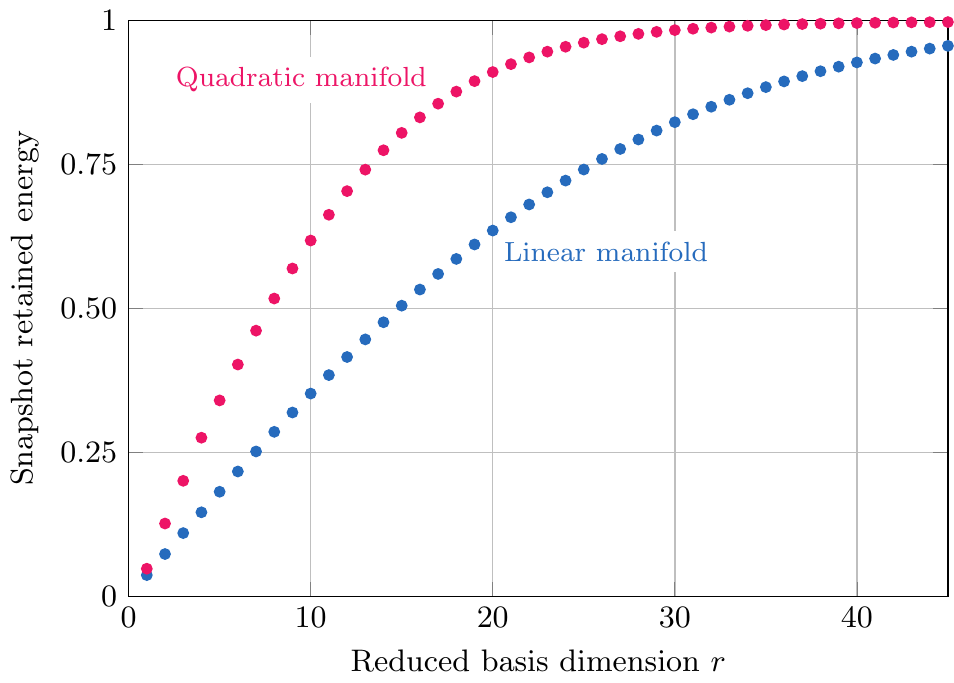}
\end{subfigure}
\caption{The snapshot retained energy spectrum of the training dataset computed with \eqref{eq:snapshot_nrg1} and \eqref{eq:snapshot_nrg2}. The regularizer to compute $\overline{\mathbf{V}}$ in \eqref{eq:least-squares_reg1} is chosen to be $\gamma=10^9$.}
\label{fig:1d_energy}
\end{figure}

Let us consider the one-dimensional linear transport equation
\begin{equation}
    \dfrac{\partial }{\partial t} s(x,t) + c\dfrac{\partial}{\partial x} s(x,t) = 0, \quad x \in \mathbb{R},
    \label{eq:1d_advection}
\end{equation}
with time $t \in (0,\infty )$, whose exact solution is given by $s(x,t)=s_0(x - c t)$, where $c$ is a constant advection velocity and $s_0(x)$ is the specified initial state at $t=0$. The state $s(x,t)$ might, for example, represent the concentration of a pollutant being advected in the one-dimensional flow with constant velocity $c$.
We consider an initial condition of the form
\begin{equation}
\label{eq:hx1}
    s_0(x) := s(x,0) = \dfrac{1}{\sqrt{0.0002 \pi}} \exp\left( -\dfrac{(x-\mu)^2}{0.0002}\right), \quad x \in \mathbb{R},
\end{equation}
where $\mu$ is a one-dimensional parameter in the parameter space $\mathcal{P}$.

In this experiment a training dataset is built by uniformly sampling the exact solution in the space-time domain $(0,1)\times(0,0.1)$ considering $n=2^{12}$ degrees of freedom for the spatial discretization and $k=2000$ time instances. The time derivatives $\partial s(x,t)/\partial t$ are computed in closed form. The parameter $\mu$ varies in the one-dimensional parameter space $\mathcal{P} = [0.05, 0.25]$. Figure~\ref{fig:1d_ic} shows the initial condition for $\mu=0.10$. We consider $\mathcal{N}_\text{test}=50$ testing-parameter instances, randomly sampled over $\mathcal{P}$, to evaluate the constructed reduced-order models on unseen data (see also \cite{fresca2021comprehensive}). We reiterate the importance of properly centering the training data about an appropriate reference state $\mathbf{s}_\text{ref}$.
The choice of reference state can affect accuracy and stability of data-driven reduced-order models based on quadratic manifolds.
Here, we center the training data about its time-averaged mean value. The normalized singular values of the centered data matrix are plotted in Figure~\ref{fig:1d_svs}. The decay of the singular values suggests that a reduced space of dimension $r \simeq 200$ is necessary for approximating the trajectories $s(x,t)$ in a linear subspace with a projection error of $10^{-8}$ in the Euclidean norm.

\begin{figure}[tbp]
\centering \small
\begin{subfigure}[t]{.49\linewidth}
\includegraphics[scale=1]{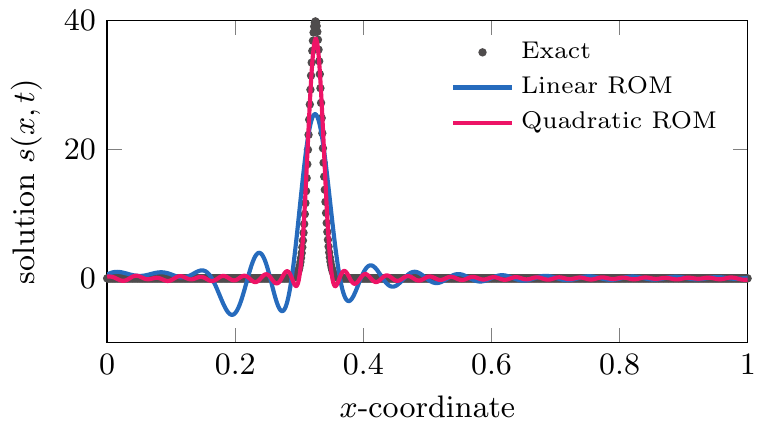}
\caption{$t=.02$}
\end{subfigure}
\begin{subfigure}[t]{.49\linewidth}
\includegraphics[scale=1]{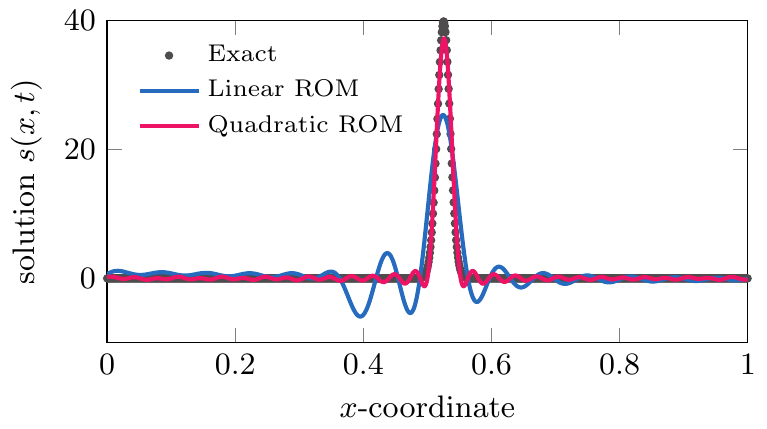}
\caption{$t=.04$}
\end{subfigure}
\begin{subfigure}[t]{.49\linewidth}
\includegraphics[scale=1]{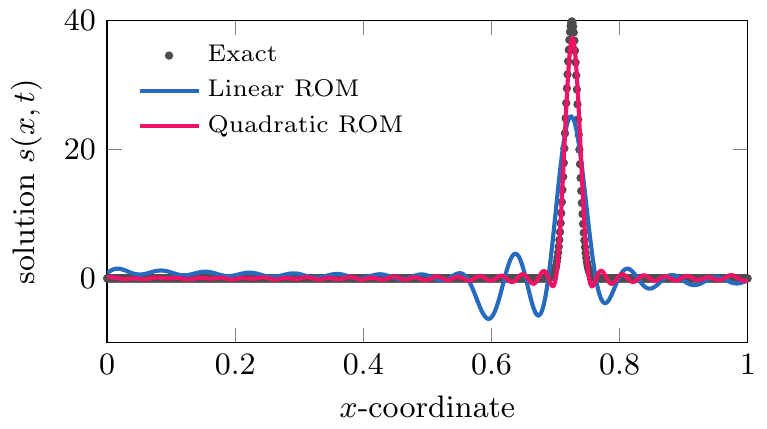}
\caption{$t=.06$}
\end{subfigure}
\begin{subfigure}[t]{.49\linewidth}
\includegraphics[scale=1]{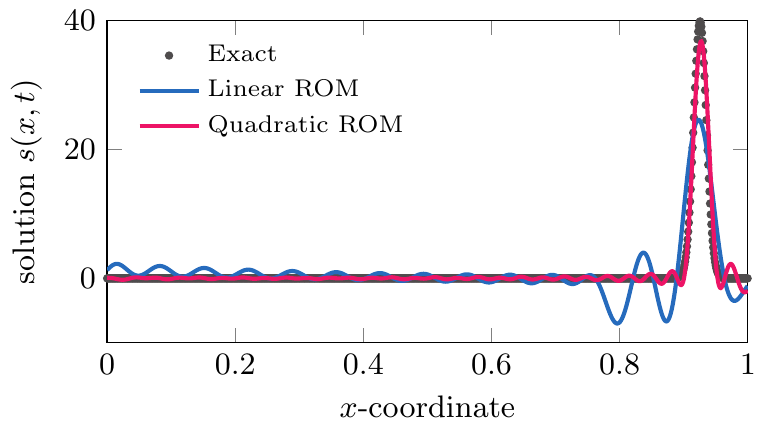}
\caption{$t=.08$}
\end{subfigure}
\caption{Comparison of the simulated solutions at four different time instances $t \in \{.02, .04, .06, .08\}$ computed by solving the linear-subspace reduced model for $r=29$ (blue lines) and quadratic-subspace reduced model with the same value of $r$, using all $r(r+1)/2$ quadratic modes (red lines) with $c=10$ at $\mu_\text{test}=0.12547$. The exact solution (black markers) is shown for reference. The plots show the reconstructed solutions for reduced models where the basis accounts for $80.9\%$ of the snapshot energy with the linear manifold approach and $98.1\%$ with the quadratic manifold approach, as computed by \eqref{eq:snapshot_nrg1} and \eqref{eq:snapshot_nrg2}, respectively. }
\label{fig:1d_2}
\end{figure}

\begin{figure}[tbp]
\centering \small
\begin{subfigure}[t]{0.28\textwidth}
\includegraphics[width=\linewidth]{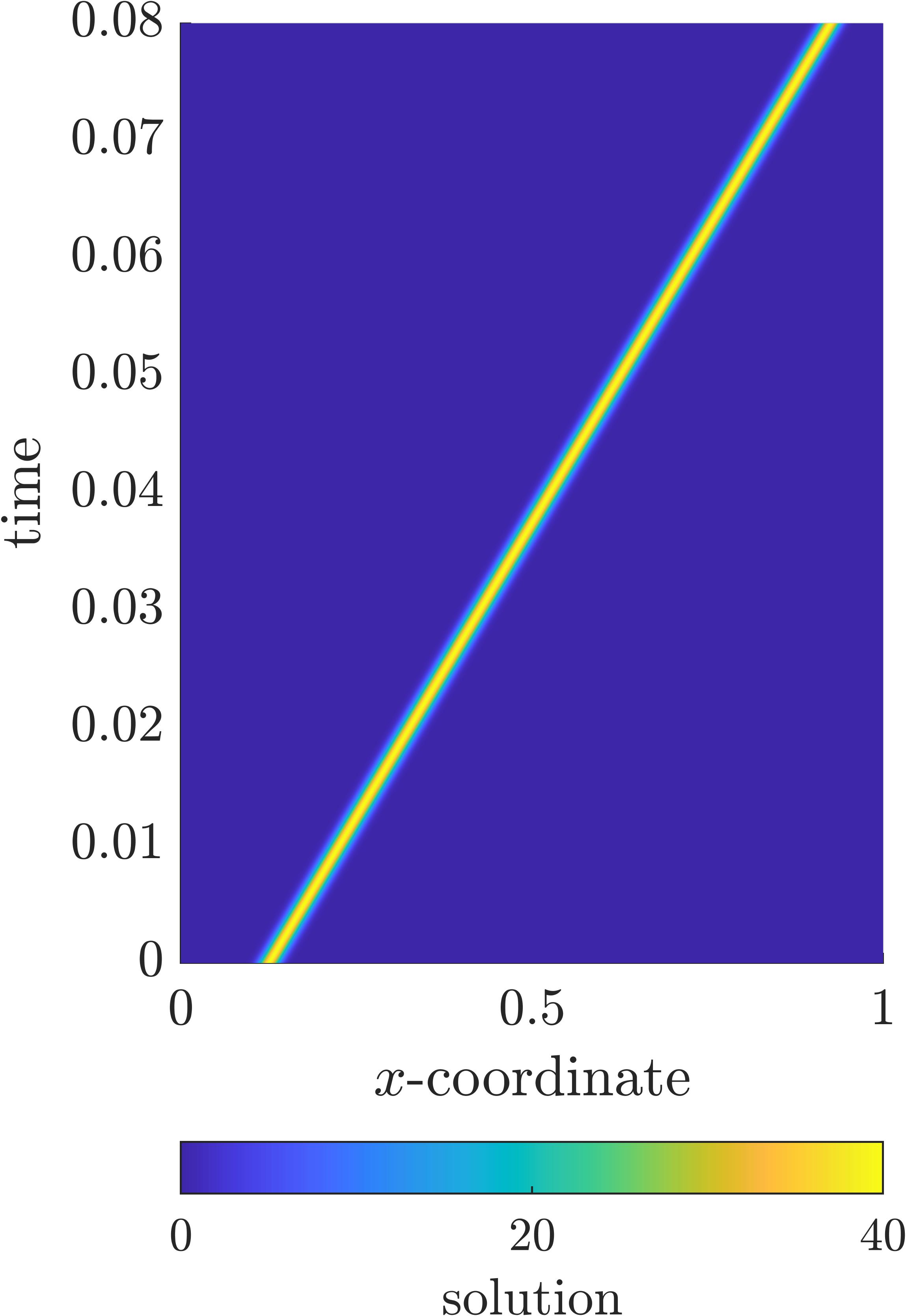}
\caption{Exact}
\end{subfigure}
\begin{subfigure}[t]{0.28\textwidth}
\includegraphics[width=\linewidth]{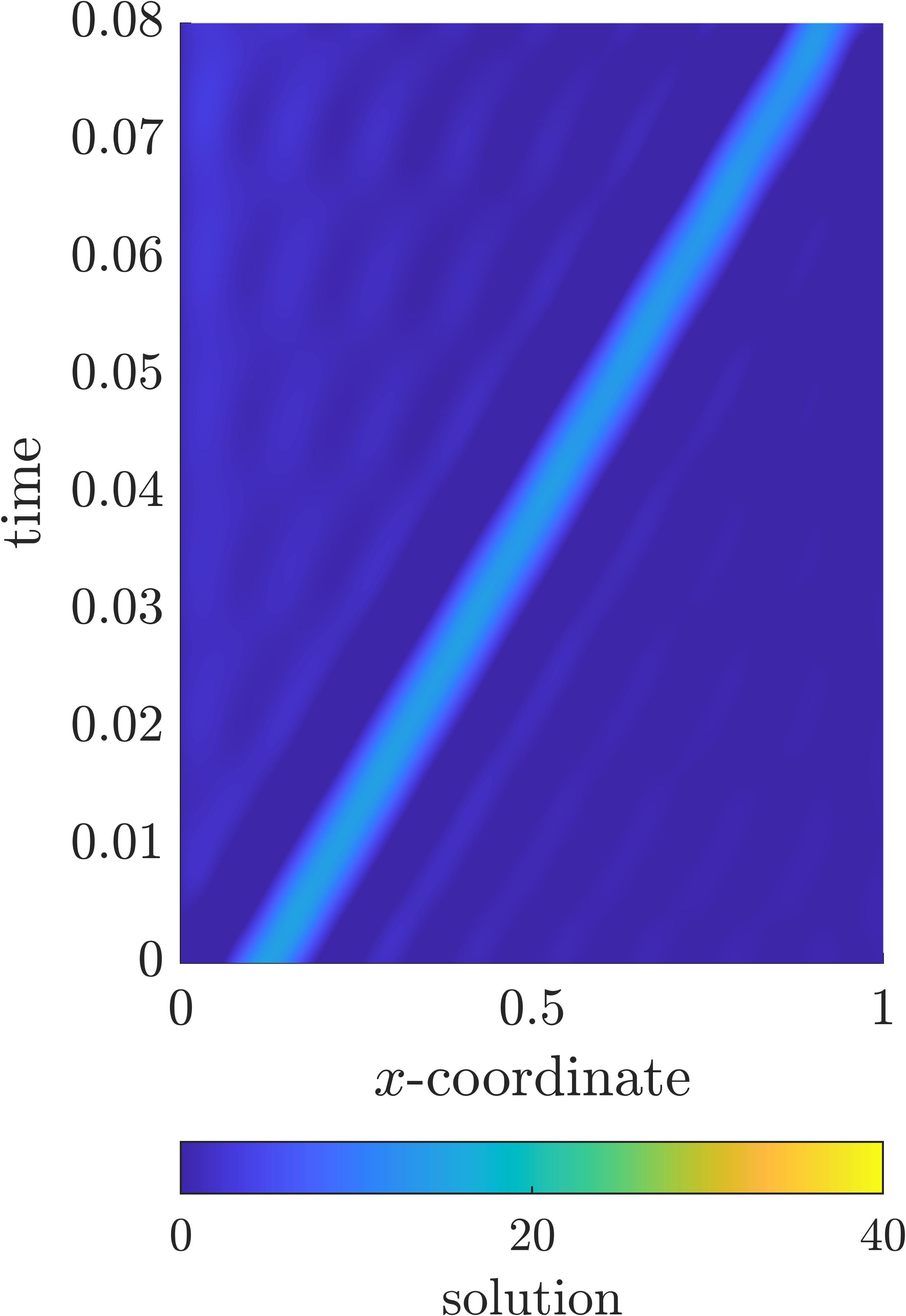}
\caption{Linear manifold}
\end{subfigure}
\begin{subfigure}[t]{0.28\textwidth}
\includegraphics[width=\linewidth]{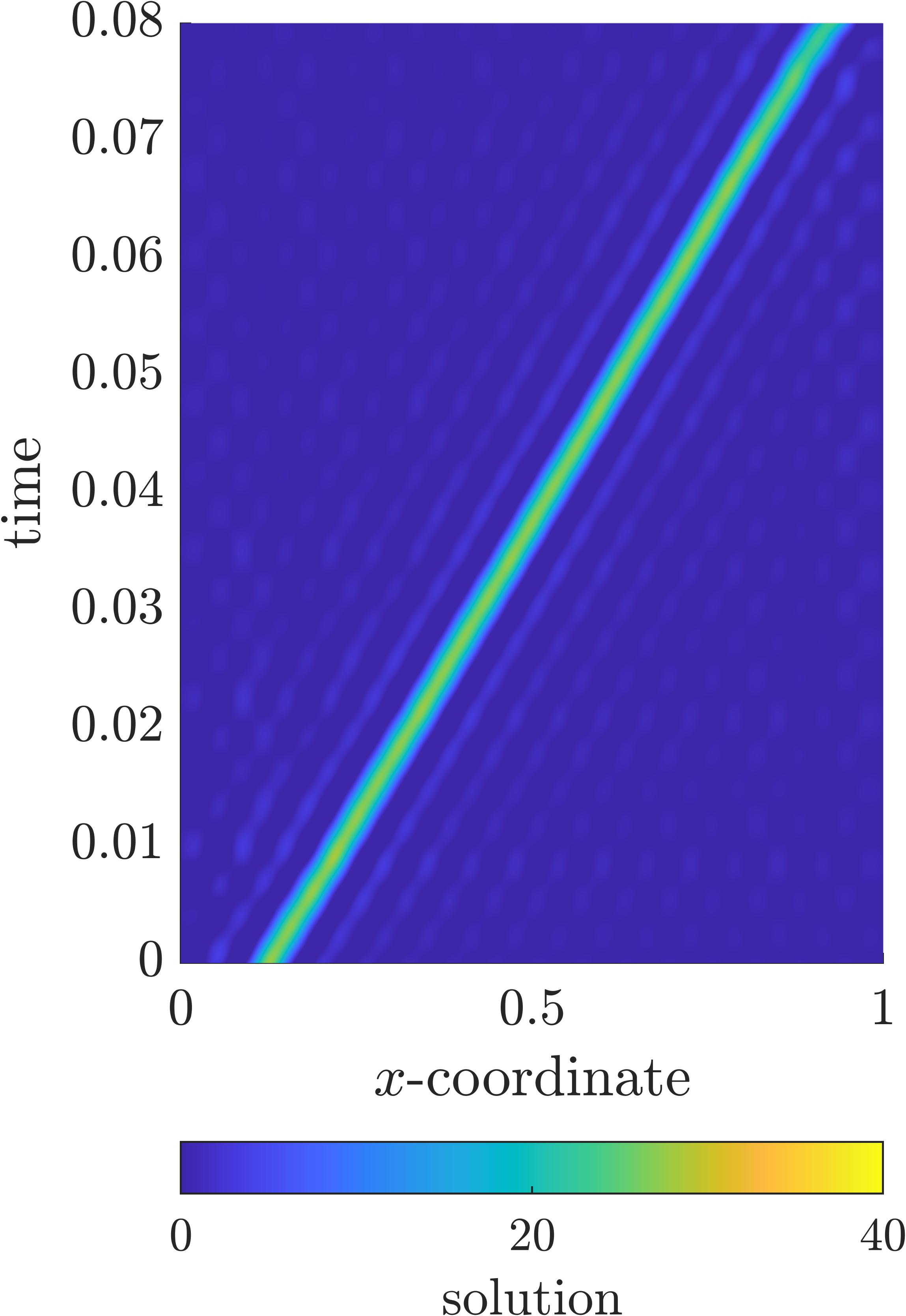}
\caption{Quadratic manifold}
\end{subfigure}
\caption{Comparison of the evolution of the solution field $s(x,t)$ through space-time at testing parameter $\mu_\text{test}=0.12547$ for (a)~the exact solution and the operator inference reduced-order models with (b)~a linear manifold and (c)~a quadratic manifold using all $r(r+1)/2$ coordinate combinations. The reduced models employ $r=15$ POD modes capturing 50.5\% and 80.5\% of the snapshot energy for the linear and quadratic manifold approaches, respectively.}
\label{fig:1d_r=15}
\end{figure}

\begin{figure}[tbp]
\centering \small
\begin{subfigure}[t]{0.28\textwidth}
\includegraphics[width=\linewidth]{./figures/fom.jpg}
\caption{Exact}
\end{subfigure}
\begin{subfigure}[t]{0.28\textwidth}
\includegraphics[width=\linewidth]{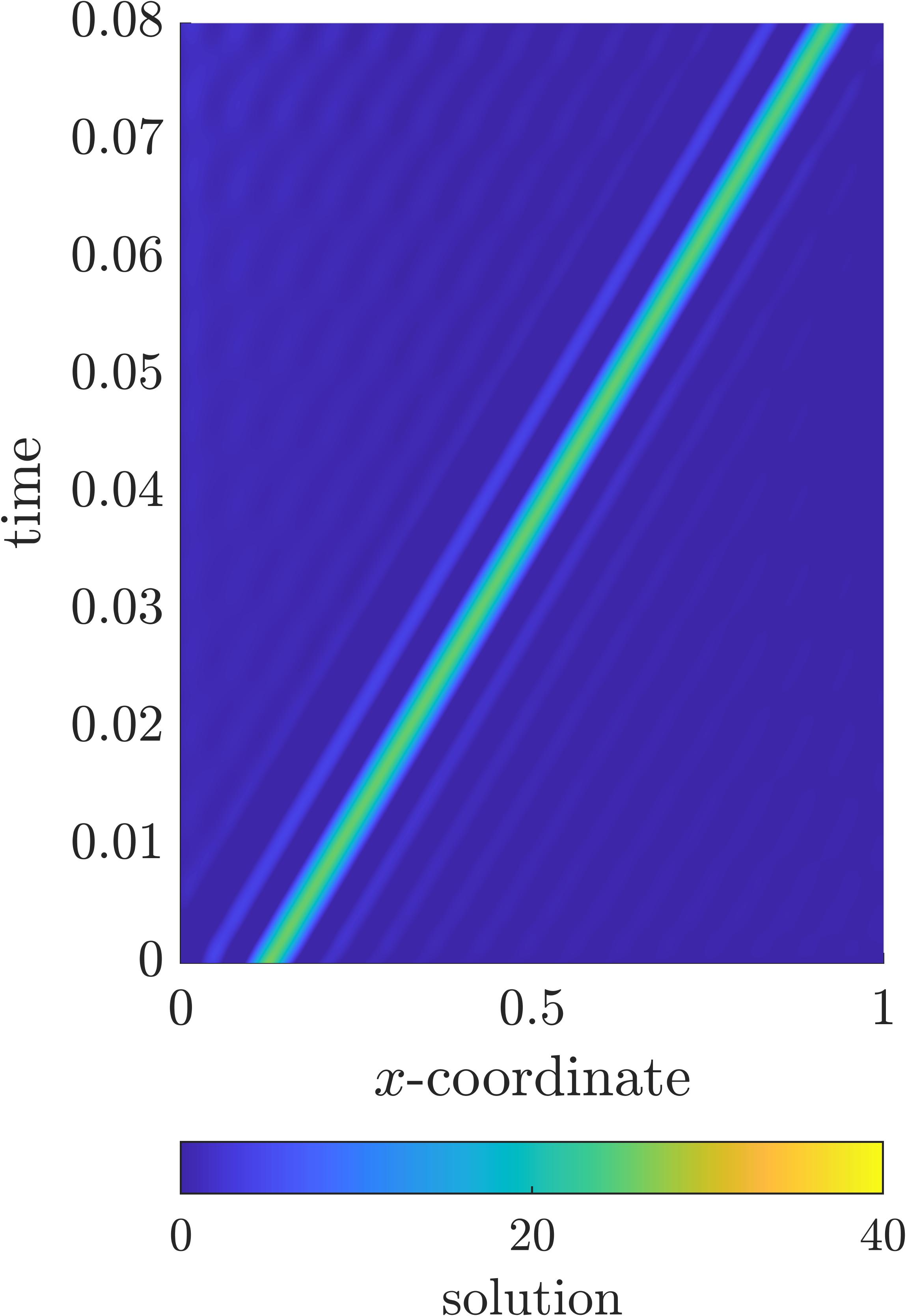}
\caption{Linear manifold}
\end{subfigure}
\begin{subfigure}[t]{0.28\textwidth}
\includegraphics[width=\linewidth]{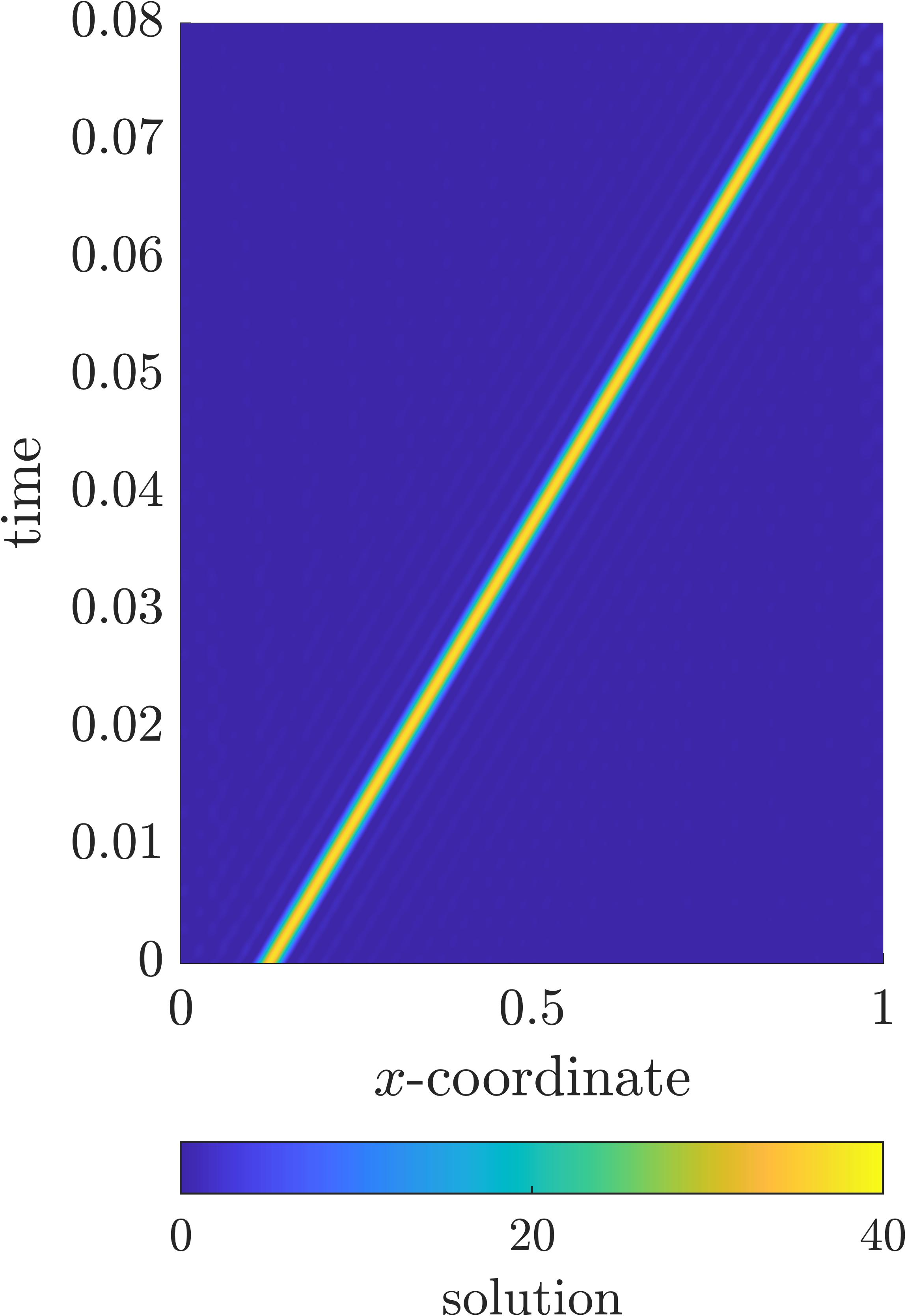}
\caption{Quadratic manifold}
\end{subfigure}
\caption{Comparison of the evolution of the solution field $s(x,t)$ through space-time at testing parameter $\mu_\text{test}=0.12547$ for (a)~the exact solution and the operator inference reduced-order models with (b)~a linear manifold and (c)~a quadratic manifold (using all $r(r+1)/2$ coordinate combinations). The reduced models employ $r= 29$ POD modes capturing 80.9\% and 98.1\% of the snapshot energy for the linear and quadratic manifold approaches, respectively.}
\label{fig:1d_r=29}
\end{figure}

A typical approach for linear dimensionality reduction techniques is to choose $r$ such that approximation in the $r$-dimensional basis $\mathbf{V}\in \mathbb{R}^{n \times r}$ yields
\begin{equation}
\dfrac{\| \mathbf{V}\mathbf{V}^\top (\mathbf{S}-\mathbf{S}_\text{ref}) \|_F^2}{\| \mathbf{S}-\mathbf{S}_\text{ref} \|_F^2} = \dfrac{\sum_{i=1}^r\sigma_i^2}{\sum_{i=1}^k\sigma_i^2} > \kappa
\label{eq:snapshot_nrg1}
\end{equation}
where $\kappa$ is a user-specified tolerance $\kappa$ and the $\sigma_i^2$ are the squared singular values of the shifted data matrix. The left-hand side of \eqref{eq:snapshot_nrg1} is often referred to as the ``relative cumulative energy'' of the system captured by $r$ POD modes. For the proposed manifold approach an equivalent snapshot ``retained energy'' metric can be devised for choosing the reduced basis dimension:
\begin{equation}
    \dfrac{\| \mathbf{V}\mathbf{V}^\top (\mathbf{S}-\mathbf{S}_\text{ref}) + \overline{\mathbf{V}} ( \mathbf{V}^\top (\mathbf{S}-\mathbf{S}_\text{ref}) \odot \mathbf{V}^\top (\mathbf{S}-\mathbf{S}_\text{ref}) ) \|_F^2}{\| \mathbf{S}-\mathbf{S}_\text{ref} \|_F^2} > \kappa,
    \label{eq:snapshot_nrg2}
\end{equation}
where $\odot$ denotes the column-wise Kronecker product of two matrices.
Figure~\ref{fig:1d_energy} plots the snapshot energy spectra of the training data computed with \eqref{eq:snapshot_nrg1} and \eqref{eq:snapshot_nrg2} as a function of the reduced basis dimension $r$. The quadratic mapping operator $\overline{\mathbf{V}}$ in \eqref{eq:snapshot_nrg2} is obtained from  \eqref{eq:least-squares_reg1} with $\gamma=10^9$, which provided adequate regularization for the reduced-order models inferred from the training data. As can be seen from Figure~\ref{fig:1d_energy}, the quadratic manifold accounts for a larger share of the total snapshot energy at a given reduced basis dimension $r$ compared to traditional linear subspace approaches.

\begin{figure}[tbp]
\centering \small
\begin{subfigure}[t]{.49\linewidth}
\hspace{-2em}\includegraphics[scale=1]{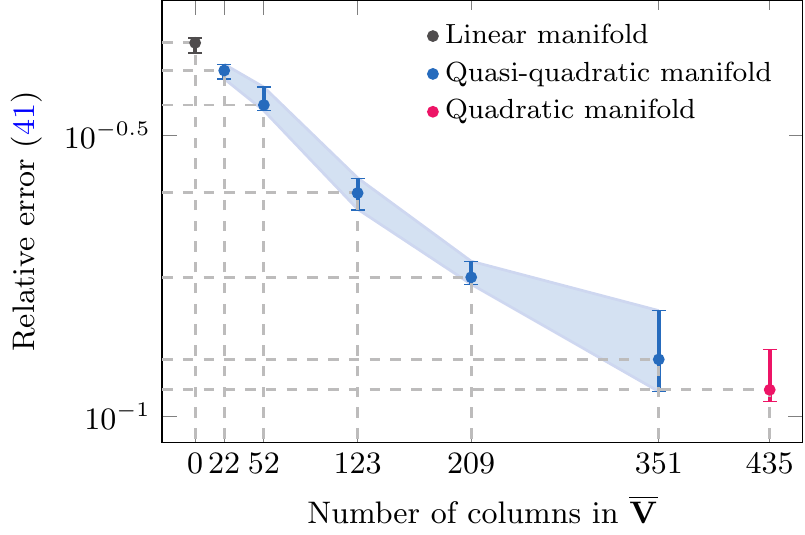}
\end{subfigure}
\caption{Relative state prediction error \eqref{eq:1d_err} as a function of the number of basis functions spanning the quadratic solution-manifold. The median and first/third quartile errors are plotted. All simulations use  a reduced-order model of dimension $r=29$.}
\label{fig:reduced_data_transport}
\end{figure}

We now employ quadratic manifolds of varying dimensions in the derivation of reduced-order models using operator inference as described in Section~\ref{subsec:opinf}.
The reduced models are inferred by solving \eqref{eq:opinf_quadratic}.
Simulations of the reduced-order models employ a semi-implicit Euler time integration scheme with time step size $\Delta t=10^{-6}$ and end time $T=0.08$.
Figure~\ref{fig:1d_2} shows representative simulation results for linear and quadratic reduced-order models with reduced dimension $r=29$, capturing 80.9\% and 98.1\% of the snapshot energy, respectively.
The exact solutions are shown for reference.
While the reconstruction for both the linear and quadratic models are subject to spurious oscillations in the solution field, they are less pronounced in the quadratic formulation.
The space-time evolution of a sample trajectory is shown in Figures~\ref{fig:1d_r=15} and \ref{fig:1d_r=29} for reduced-order models of dimension $r=15$ and $r=29$, respectively.
The performance of the operator inference quadratic reduced-order models is good in both cases, although the $r=29$ case is clearly more accurate as the approximate solution becomes visually indistinguishable from the exact solution.

We now turn to the quasi-quadratic formulation from Section \ref{subsec:quasi-quadratic}.
Recall that the proposed quadratic manifold approach does \textit{not} introduce extra degrees of freedom in the associated state dimension of the reduced-order models.
Instead, we rely on a set of basis vectors, the columns of $\overline{\mathbf{V}}$, to define the nonlinear aspects of the solution-manifold. In the baseline quadratic manifold approach (outlined in Section~\ref{subsec:nonlinear}), the number of additional basis vectors is $120$ and $435$ for the values $r=15$ and $r=29$, respectively.
By performing column selection via the regularized problem \eqref{eq:least-squares_reg2}, we seek parametrizations that use only a subset of these additional basis vectors. We can control the number of selected columns via the choice of regularization parameter $\gamma$. Algorithm~\ref{alg:sparsa} is used to solve problem \eqref{eq:least-squares_reg2} for each  value of $\gamma$.
In Figure~\ref{fig:reduced_data_transport} we plot the error in the $r=29$ reduced-order model solution as a function of the number of columns in $\overline{\mathbf{V}}$, denoted by $q$. We employ a relative error measure
\begin{equation}
    e(\mu_\text{test})= \left.\left( \sqrt{\displaystyle \sum_{j=1}^{\mathcal{N}_\text{t}} \| \mathbf{s}_j(\mu_\text{test}) - \mathbf{s}_{j,\text{approx}}(\mu_\text{test}) \|^2 } \right) \middle/ \left(\sqrt{\displaystyle \sum_{j=1}^{\mathcal{N}_\text{t}} \| \mathbf{s}_j(\mu_\text{test}) \|^2 }\right)\right. ,
    \label{eq:1d_err}
\end{equation}
where $\mathbf{s}_j(\mu_\text{test})$ denotes the exact solution sampled on the spatial discretization at time step $t_j$ using the initial condition with parameter $\mu_\text{test}$, $\mathbf{s}_{j,\text{approx}}(\mu_\text{test})$ denotes the corresponding solution predicted with the reduced model, and
$\mathcal{N}_\text{t}$ is the number of time instances at which the solution is computed for each testing parameter.\footnote{While this indicator is frequently used to assess the performance of reduced-order models, see \cite{fresca2021comprehensive} for instance, we note that it can be sensitive to misalignments of the predicted solution relative to the ground truth in transport-dominated problems.}
When $q$ is relatively small, we obtain predictions that improve only marginally over the linear parametrization.
However, in the limit as fewer of the columns are eliminated, we recover the accuracy metrics from the baseline quadratic manifold approach.
We note that for the largest value of $q$, here $q=351$, at some of the testing instances the reduced-order models were unstable and were therefore not included in error measure \eqref{eq:1d_err} and Figure~\ref{fig:reduced_data_transport}.
We emphasize that all reduced models in  Figure~\ref{fig:reduced_data_transport} have dimension $r=29$; reducing the number of columns in $\overline{\mathbf{V}}$ reduces the dimension of the reduced-order quadratic operator $\widehat{\mathbf{H}}$ and reduces the number of basis vectors that must be stored in order to reconstruct full-domain solutions, but has little impact on the cost of simulating the reduced model.

\begin{figure}[tbp]
\centering \small
\begin{subfigure}[t]{.6\linewidth}
\includegraphics[scale=1]{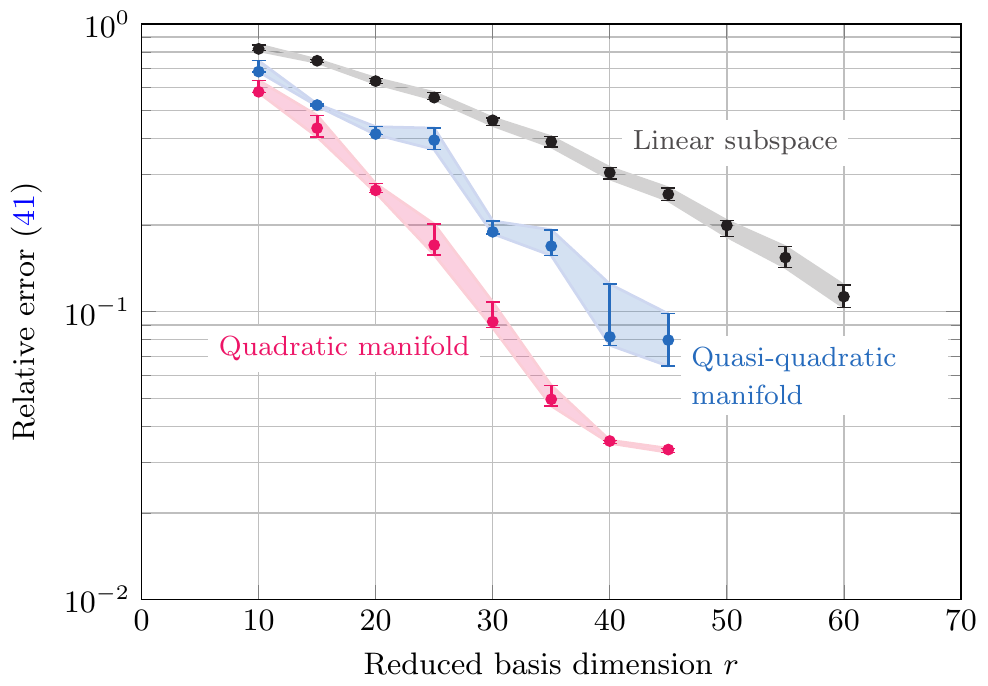}
\end{subfigure}
\caption{Relative state prediction error \eqref{eq:1d_err} as a function of the reduced basis dimension, $r$, for the testing dataset. Median and first/third quartile errors are shown. For the quasi-quadratic model, about 25\% of the $r(r+1)/2$ modes were selected.  
}
\label{fig:1d_err}
\end{figure}

We now study the behavior of the error \eqref{eq:1d_err} with respect to the reduced basis dimension $r$ across all testing parameters. Figure~\ref{fig:1d_err} shows results for the cases of linear, quadratic, and quasi-quadratic manifolds. In the quasi-quadratic case, a regularizer $\gamma$ was sought which reduces the number of columns in $\overline{\mathbf{V}}$ to the most significant 25\% of the directions in $\mathbf{V}_\perp$, that is, $q \approx r(r+1)/8$.
The plots show that increasing the reduced basis dimension leads to a lower error for all reduction methods.
Operator inference reduced-order models constructed with the conventional linear subspace approaches produce the least accurate predictions, while introducing the quadratic manifold leads to increased predictive accuracy  across the range of $r$. The quadratic manifold introduces a more rapid decay in error with $r$.  This means that the same level of accuracy can be achieved with smaller $r$. For instance, a quadratic manifold reduced model with $r=30$ achieves roughly the same level of accuracy as a linear subspace reduced model with $r=60$.
The flattening of the quadratic manifold error curve for $r>35$ is due to the conditioning of the data matrix in the operator inference regression problem.
Poor conditioning is known to introduce numerical errors at components corresponding to the less important POD basis vectors \cite{PEHERSTORFER2016196}.
Errors for the quasi-quadratic manifold approach are intermediate between the linear and quadratic manifold cases, since for this example, the reduction in the number of columns of $\overline{\mathbf{V}}$ incurs increased error in comparison to the fully quadratic representation. This is an indication that for this problem \emph{all} quadratic terms are important. In other applications, we might expect the quasi-quadratic representation with column selection to introduce efficiencies into the representation without compromising error.

\subsection{Two-dimensional wave equation}

\begin{figure}[btp]
\centering \small
\begin{subfigure}[c]{.49\linewidth}
\includegraphics[scale=1]{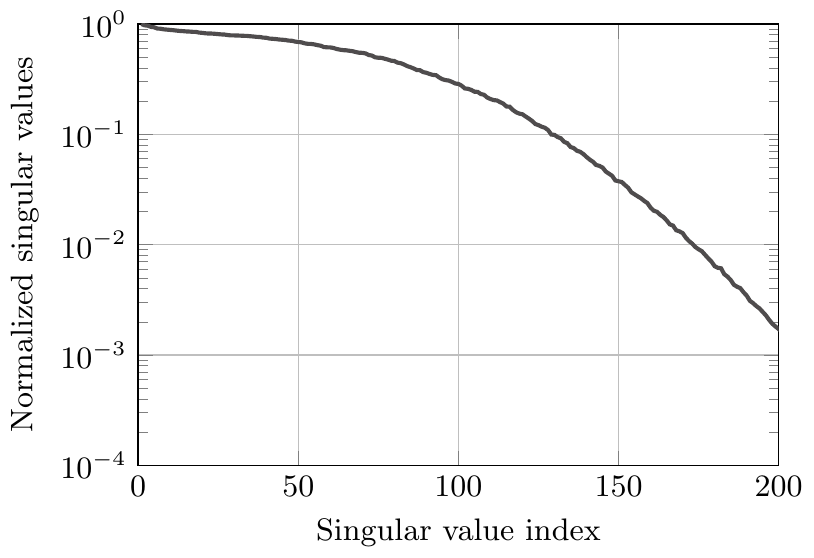}
\caption{Singular values}
\label{fig:2d_svs}
\end{subfigure}
\begin{subfigure}[c]{.49\linewidth}
\includegraphics[scale=1]{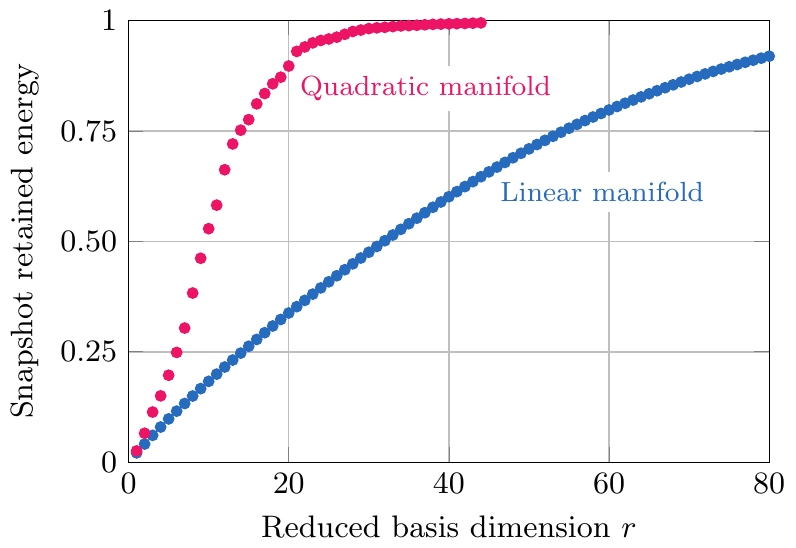}
\caption{Snapshot retained energy}
\label{fig:2d_energy}
\end{subfigure}
\caption{Two-dimensional wave equation example. Left: snapshot singular value decay. Right: snapshot retained energy spectrum of the training dataset computed with \eqref{eq:snapshot_nrg1} and \eqref{eq:snapshot_nrg2}. The regularizer to compute $\overline{\mathbf{V}}$ in \eqref{eq:least-squares_reg1} is chosen to be $\gamma=10^{-1}$.}
\end{figure}

For a second example, we consider the two-dimensional wave equation, a second-order hyperbolic linear PDE, in a rectangular domain $\Omega = [0,4\pi] \times [0,2\pi]$. The wave equation has been studied previously as a benchmark problem for nonlinear model reduction techniques in, for instance, \cite{SALVADOR20211, SHARMA2022133122}. The governing equation is given by
\begin{equation}
    \dfrac{\partial^2}{\partial t^2}s(\mathbf{x},t) = \Delta s;  \quad (\mathbf{x},t) \in \Omega \times (0,T],
    \label{eq:wave}
\end{equation}
with initial conditions
\begin{equation}
    s_0(\mathbf{x}) := s(\mathbf{x},0) = \exp \left( \dfrac{- \|\mathbf{x}-\mathbf{x}_0\|^2}{0.0072} \right); \quad \dot{s}_0(\mathbf{x}) := \dfrac{\partial}{\partial t}s(\mathbf{x},0) = 0, \quad t \in [0,T]
    \label{eq:wave_init}
\end{equation}
and homogeneous Neumann boundary conditions assigned over the entire boundary of the domain:
\begin{equation}
\nabla s(\mathbf{x}) \cdot \mathbf{n}(\mathbf{x}) = 0; \quad \mathbf{x} \in \partial \Omega,
    \label{eq:wave_bc}
\end{equation}
where $\mathbf{n}(\mathbf{x})$ is the unit outward normal vector to $\partial \Omega$. 
Here the state $s(x,t)$ might, for example, represent a physical displacement field.
The simulations are performed up to a final time $T=10$. The initial condition, a Gaussian pulse, will trigger the propagation of spherical waves about a point $\mathbf{x}_0$. Spatial discretization is performed using the finite element method on a uniform mesh with quadratic triangular elements, resulting in snapshots with $n=321,201$ entries each. For discretization in time we use Newmark implicit time integration in generating $k=1000$ training snapshots. The training snapshots are centered using their time-averaged mean value. Figure~\ref{fig:2d_svs} plots the singular values of the centered snapshot matrix. The singular values drop by about three orders in magnitude at a reduced basis dimension of $r=200$. The slow decay of the singular values is a reflection of the transport-dominated dynamics of the problem.

We apply the proposed quadratic manifold operator inference approach from Section \ref{sec:model_reduction} to this problem.
The quadratic mapping operator, $\overline{\mathbf{V}}$, is obtained through solving \eqref{eq:least-squares_reg1}. The numerical results presented here are for a regularization hyperparameter of $\gamma=10^{-1}$ in \eqref{eq:least-squares_reg1}.
This value was selected by minimizing the relative state error, \eqref{eq:1d_err}, over the training data for the case of $r=40$.
The snapshot retained energy, as defined in \eqref{eq:snapshot_nrg1} and \eqref{eq:snapshot_nrg2}, is plotted in Figure \ref{fig:2d_energy}.
As in the previous example, for a given reduced basis dimension $r$ we obtain a clear increase in the representation power of the quadratic manifold over the linear POD subspace.
For this problem, the quadratic manifold is computed for reduced basis dimensions $r\leq 44$. Beyond $r=44$, at which point the retained energy is close to 100\%, the least-squares problem associated with inferring $\overline{\mathbf{V}}$ becomes under-determined.

\begin{figure}[tbp]
\small
\begin{subfigure}[b]{0.14\textwidth}
\caption*{}
\end{subfigure}
\begin{subfigure}[b]{0.84\textwidth}
\centering
\includegraphics[width=.45\linewidth]{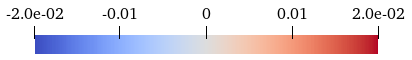}
\end{subfigure} \\[0.5em]
\begin{subfigure}[b]{0.14\textwidth}
\captionsetup{justification=raggedright, singlelinecheck=false}
\caption*{\scriptsize Reference}
\end{subfigure}
\begin{subfigure}[b]{0.28\textwidth}
\centering
\includegraphics[width=.98\linewidth]{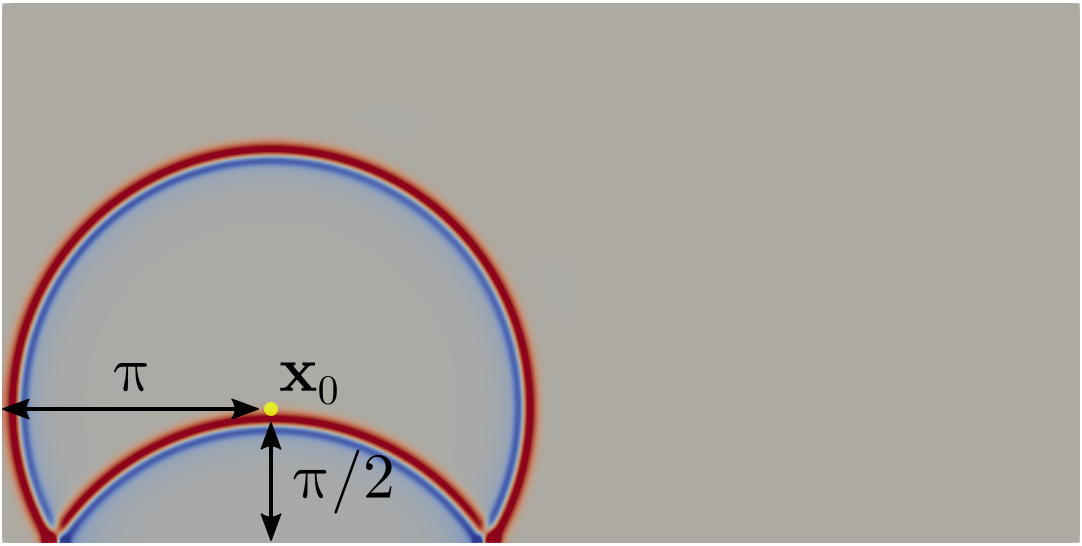}
\end{subfigure}
\begin{subfigure}[b]{0.28\textwidth}
\centering
\includegraphics[width=.98\linewidth]{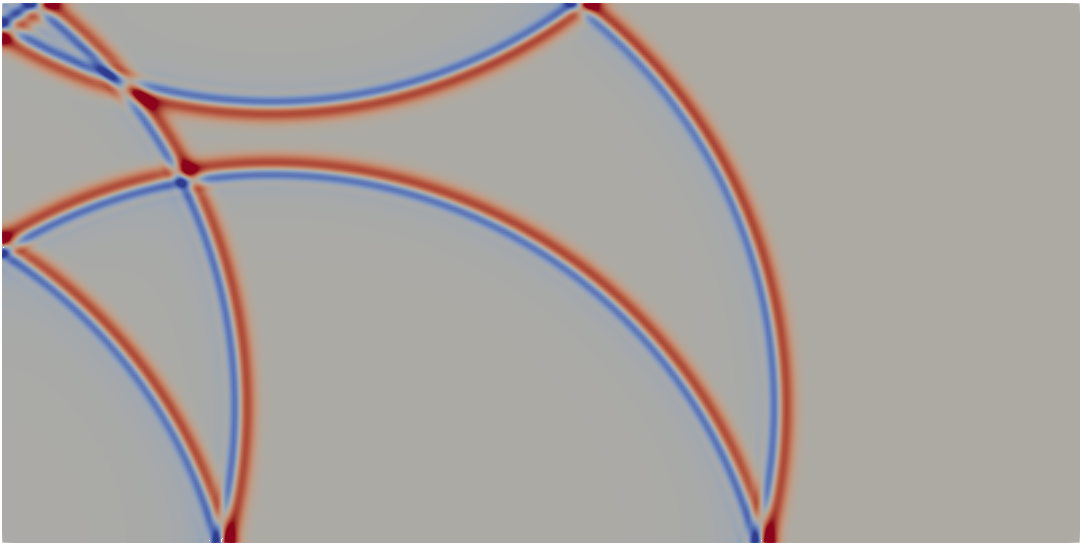}
\end{subfigure}
\begin{subfigure}[b]{0.28\textwidth}
\centering
\includegraphics[width=.98\linewidth]{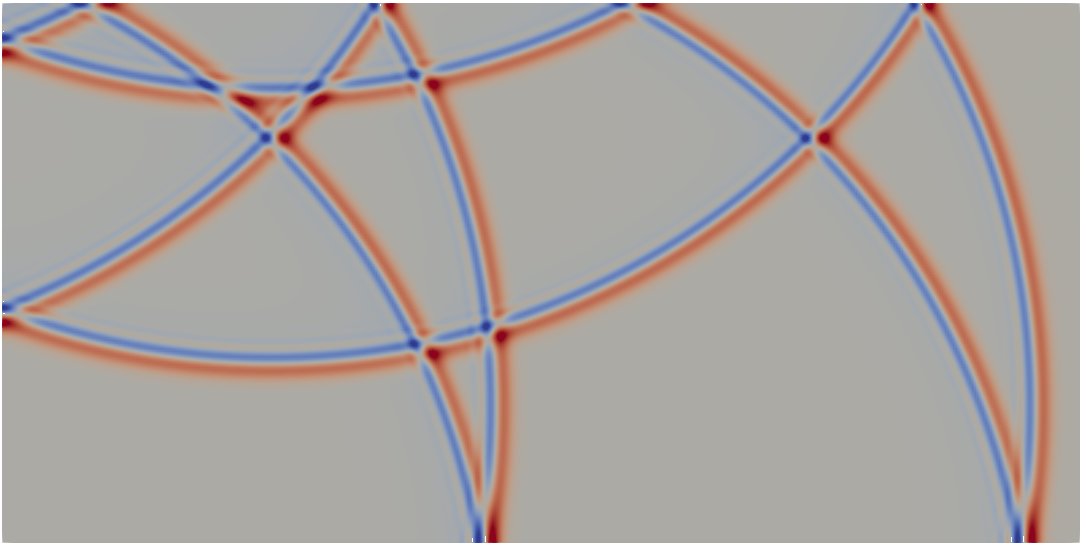}
\end{subfigure} \\[0.5em]
\begin{subfigure}[b]{0.14\textwidth}
\captionsetup{justification=raggedright, singlelinecheck=false}
\caption*{\scriptsize Operator Inference ROM (linear)}
\end{subfigure}
\begin{subfigure}[b]{0.28\textwidth}
\centering
\includegraphics[width=.98\linewidth]{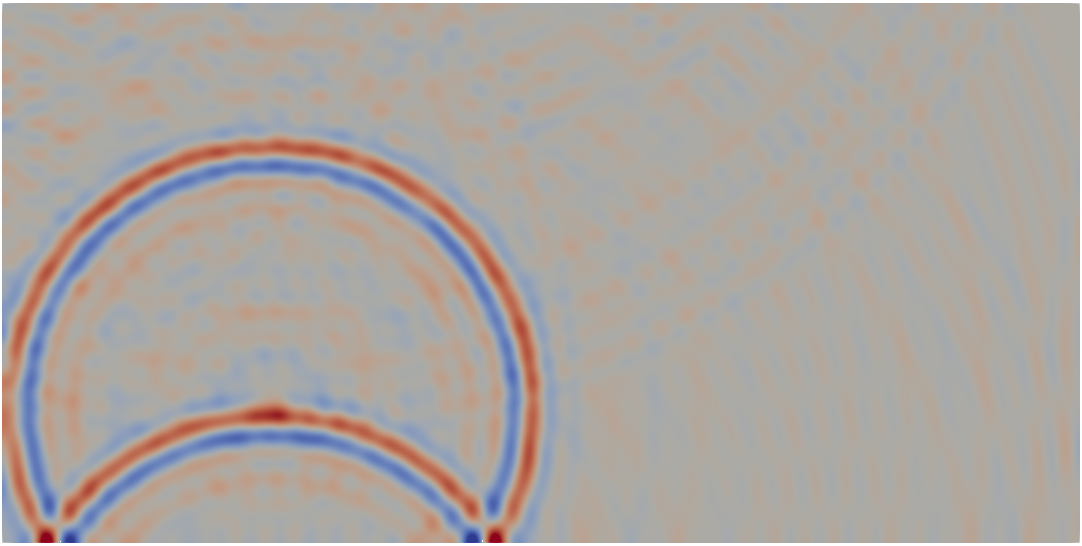}
\end{subfigure}
\begin{subfigure}[b]{0.28\textwidth}
\centering
\includegraphics[width=.98\linewidth]{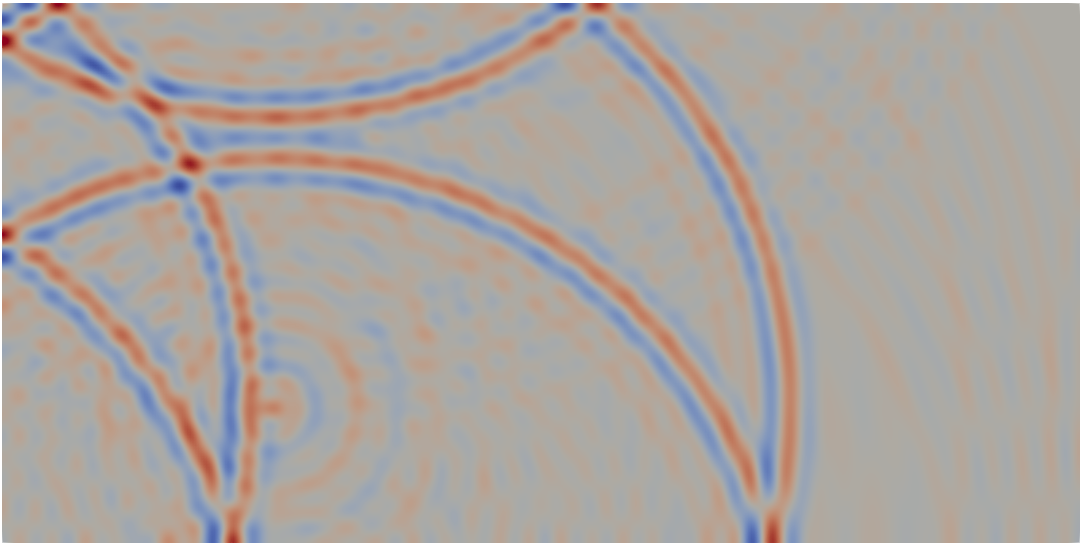}
\end{subfigure}
\begin{subfigure}[b]{0.28\textwidth}
\centering
\includegraphics[width=.98\linewidth]{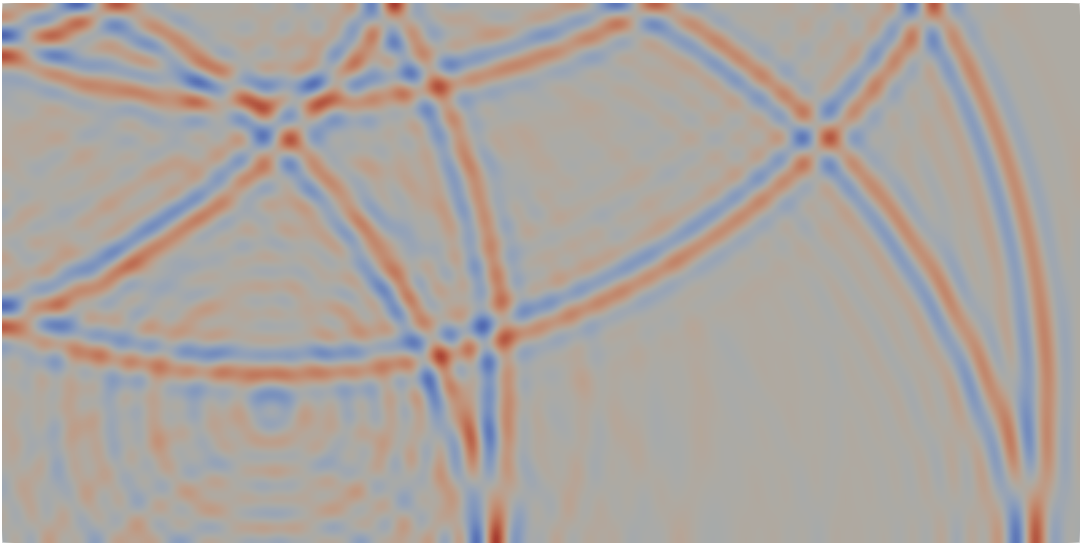}
\end{subfigure} \\[0.5em]
\begin{subfigure}[b]{0.14\textwidth}
\captionsetup{justification=raggedright, singlelinecheck=false}
\caption*{\scriptsize Operator Inference ROM (quadratic)}
\end{subfigure}
\begin{subfigure}[b]{0.28\textwidth}
\centering
\includegraphics[width=.98\linewidth]{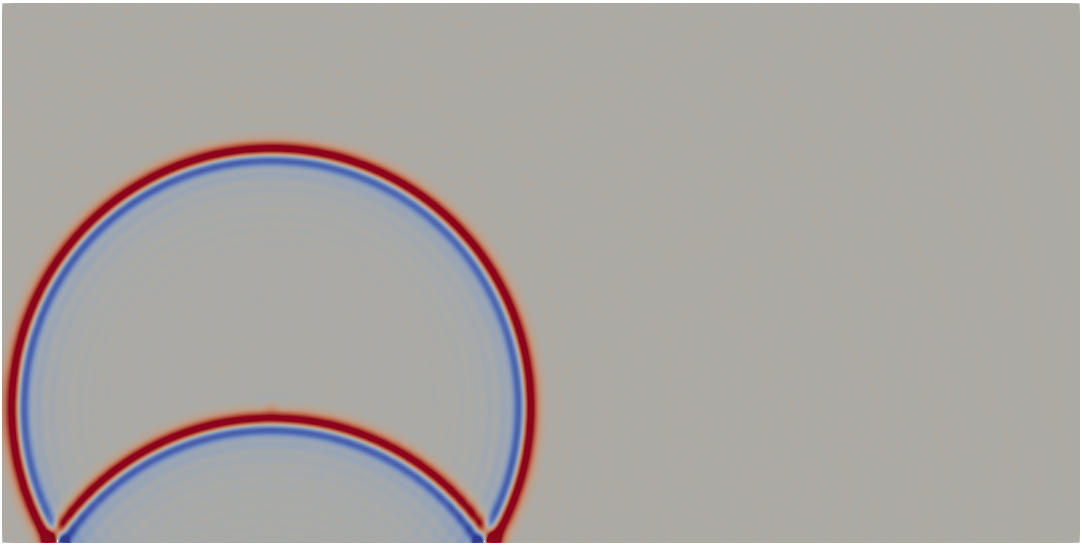}
\end{subfigure}
\begin{subfigure}[b]{0.28\textwidth}
\centering
\includegraphics[width=.98\linewidth]{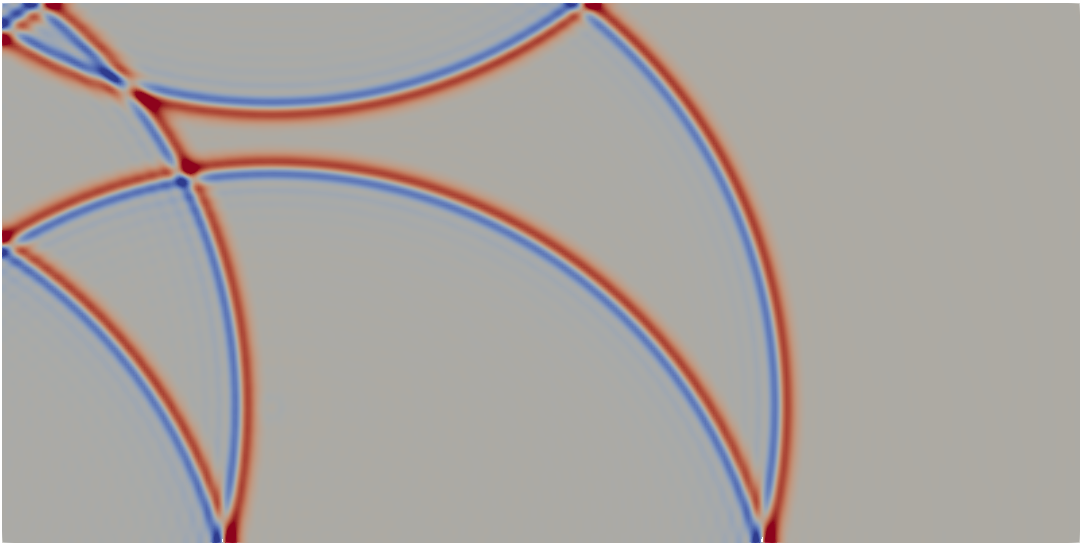}
\end{subfigure}
\begin{subfigure}[b]{0.28\textwidth}
\centering
\includegraphics[width=.98\linewidth]{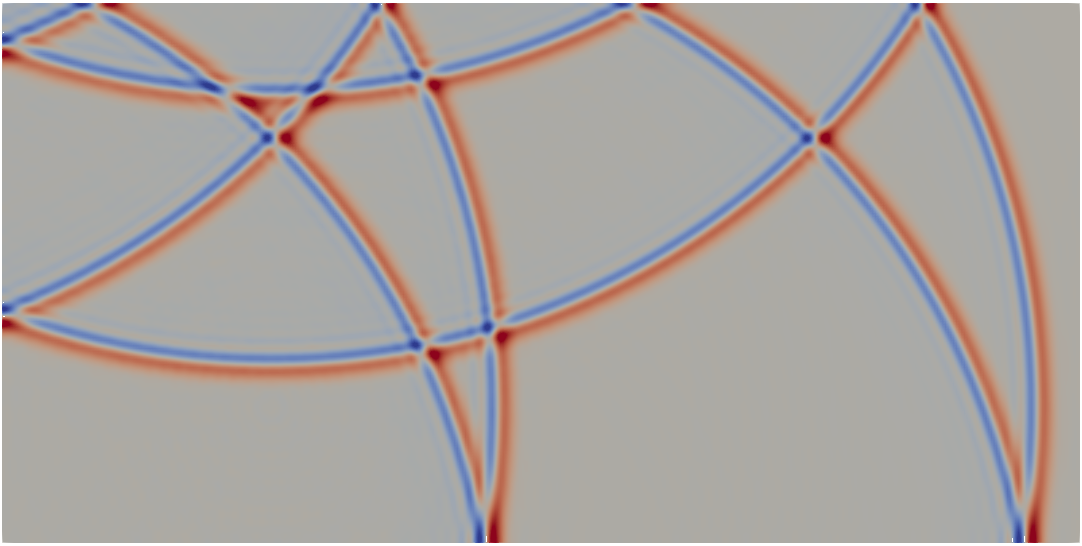}
\end{subfigure} \\[0.5em]
\begin{subfigure}[b]{0.14\textwidth}
\caption*{}
\end{subfigure}
\begin{subfigure}[b]{0.28\textwidth}
\caption{$t=3$}
\end{subfigure}
\begin{subfigure}[b]{0.28\textwidth}
\caption{$t=6$}
\end{subfigure}
\begin{subfigure}[b]{0.28\textwidth}
\centering
\caption{$t=9$}
\end{subfigure}
\caption{Comparison of the evolution of the solution field $s(x,t)$ at selected time steps $t \in \{3,6,9\}$ for the full-order model (top row) and the non-intrusive operator inference reduced-order models with approximation in a linear subspace (middle row) and the proposed quadratic manifold formulation (bottom row). Both reduced model simulations are carried out using $r=40$ POD modes, accounting for 60.1\% and 98.1\% of the snapshot retained energy, respectively, as computed by \eqref{eq:snapshot_nrg1}, \eqref{eq:snapshot_nrg2}.}
\label{fig:2d}
\end{figure}

\begin{figure}[tbp]
\centering \small
\begin{subfigure}[t]{.75\linewidth}
\includegraphics[scale=1]{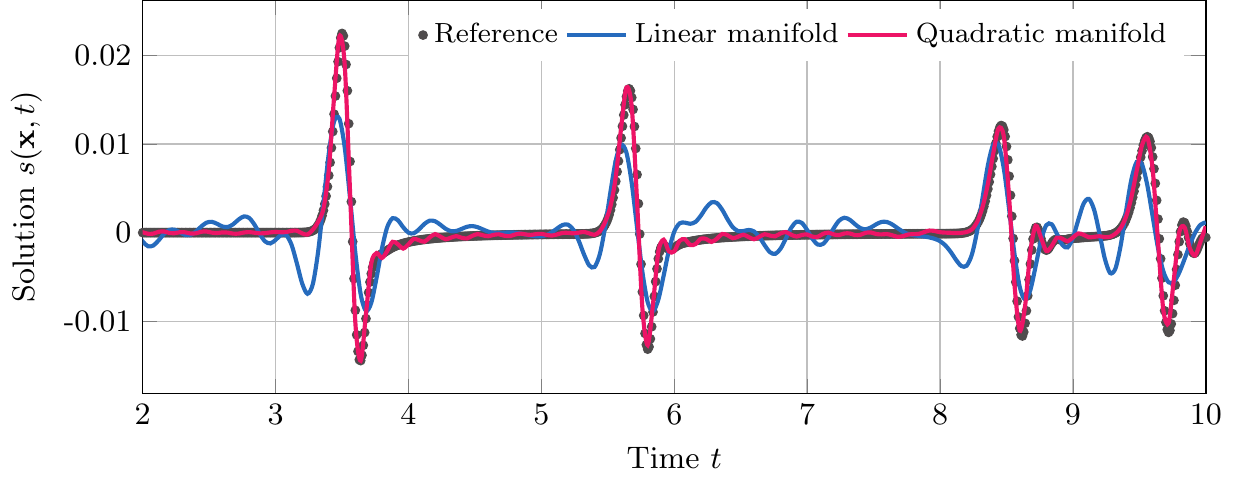}
\end{subfigure}
\caption{Traces of the solution $s(\mathbf{x},t)$ through time monitored at the center of the computational domain for the linear and quadratic operator inference models. The trace of the full-order model (black dots) with $n=321,201$ is shown for reference. The reduced basis dimension for both reduced-order models is $r=40$.}
\label{fig:trace}
\end{figure}

\begin{figure}[tbp]
\begin{subfigure}[b]{0.14\textwidth}
\caption*{}
\end{subfigure}
\begin{subfigure}[b]{0.84\textwidth}
\centering
\includegraphics[width=.45\linewidth]{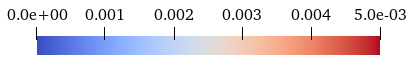}
\end{subfigure} \\[0.5em]
\begin{subfigure}[b]{0.14\textwidth}
\captionsetup{justification=raggedright, singlelinecheck=false}
\caption*{\scriptsize Absolute error: Operator Inference ROM (linear)}
\end{subfigure}
\begin{subfigure}[b]{0.28\textwidth}
\centering
\includegraphics[width=.98\linewidth]{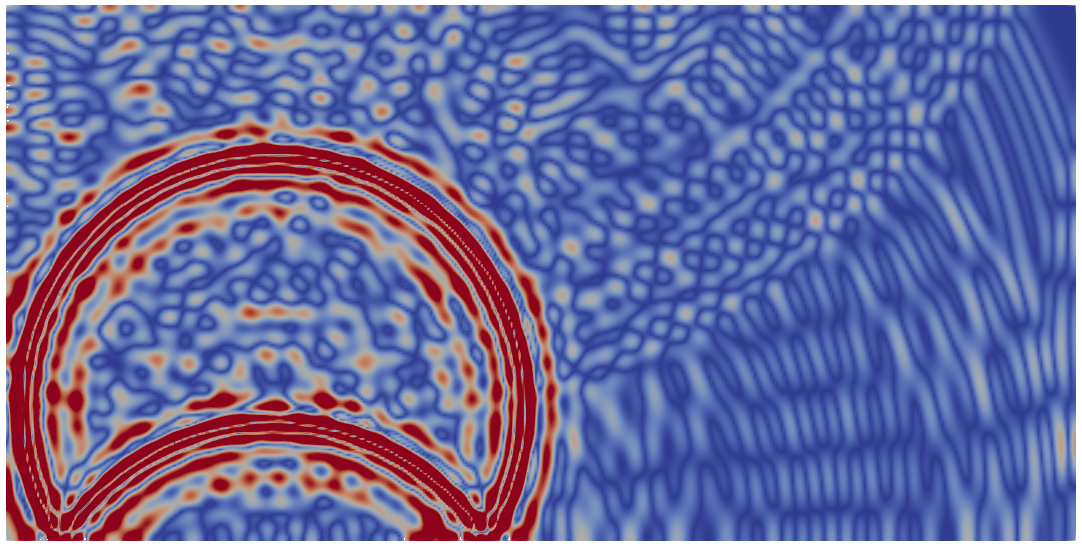}
\end{subfigure}
\begin{subfigure}[b]{0.28\textwidth}
\centering
\includegraphics[width=.98\linewidth]{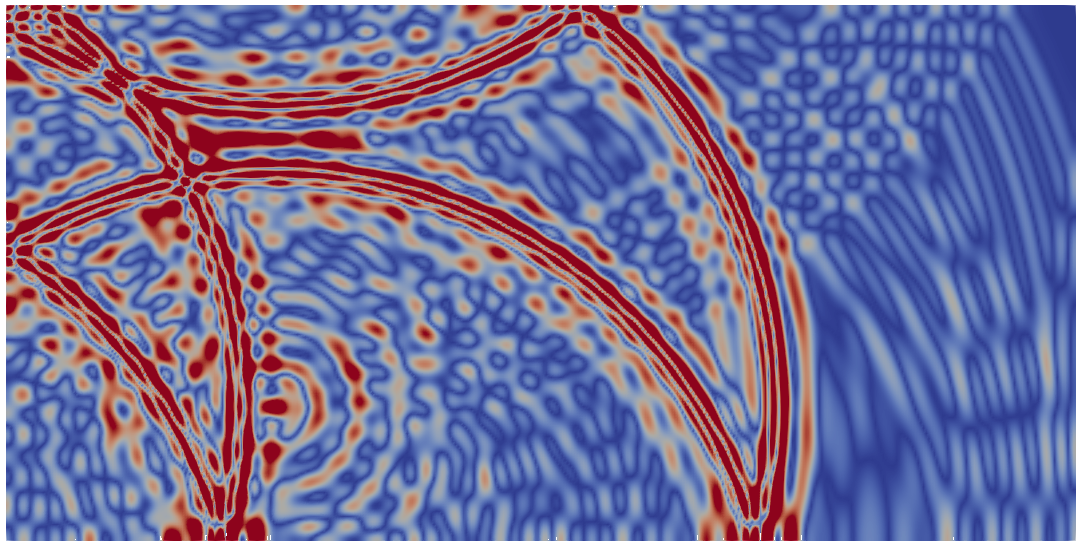}
\end{subfigure}
\begin{subfigure}[b]{0.28\textwidth}
\centering
\includegraphics[width=.98\linewidth]{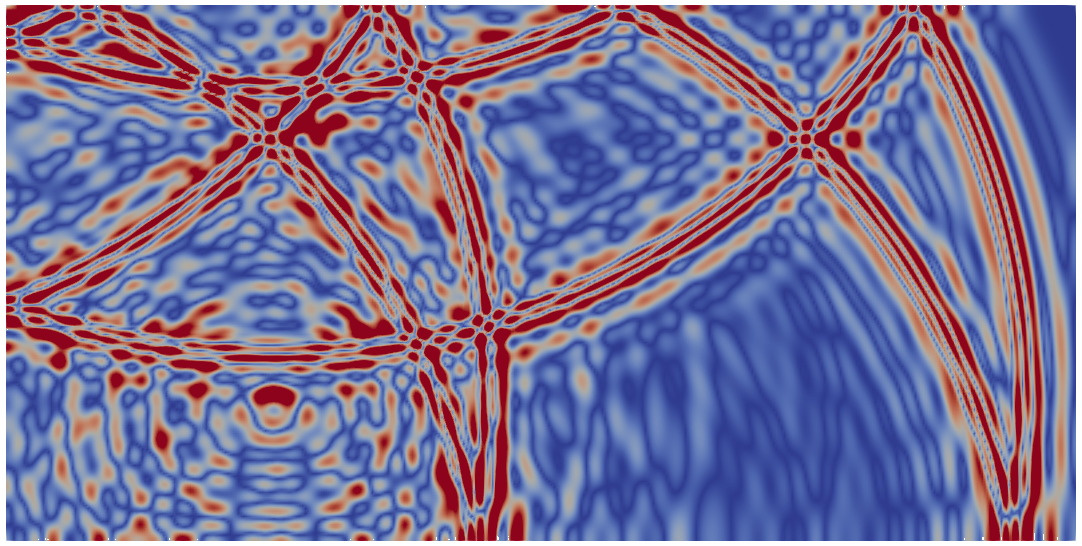}
\end{subfigure} \\[0.5em]
\begin{subfigure}[b]{0.14\textwidth}
\captionsetup{justification=raggedright, singlelinecheck=false}
\caption*{\scriptsize Absolute error: Operator Inference ROM (quadratic)}
\end{subfigure}
\begin{subfigure}[b]{0.28\textwidth}
\centering
\includegraphics[width=.98\linewidth]{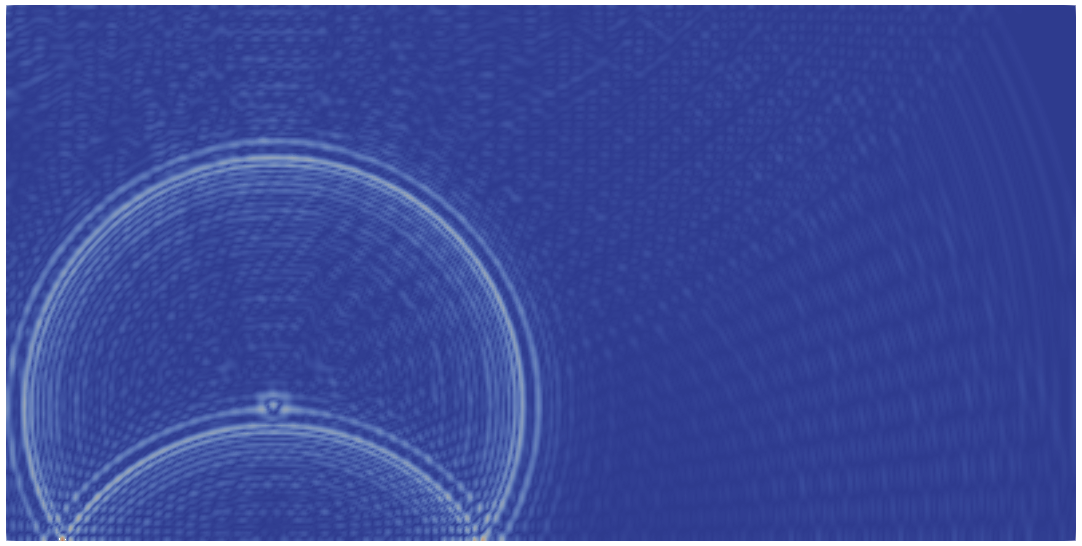}
\end{subfigure}
\begin{subfigure}[b]{0.28\textwidth}
\centering
\includegraphics[width=.98\linewidth]{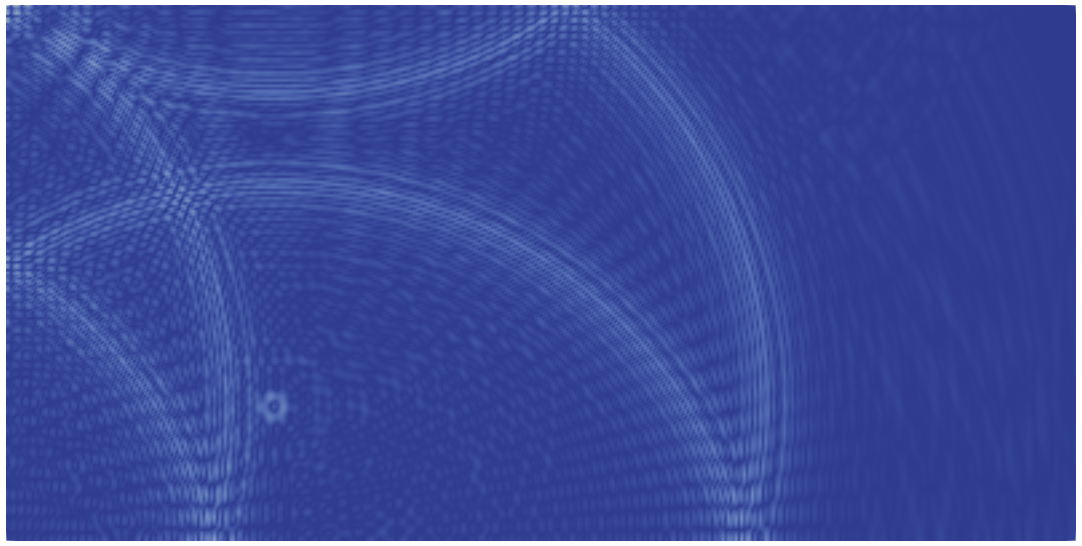}
\end{subfigure}
\begin{subfigure}[b]{0.28\textwidth}
\centering
\includegraphics[width=.98\linewidth]{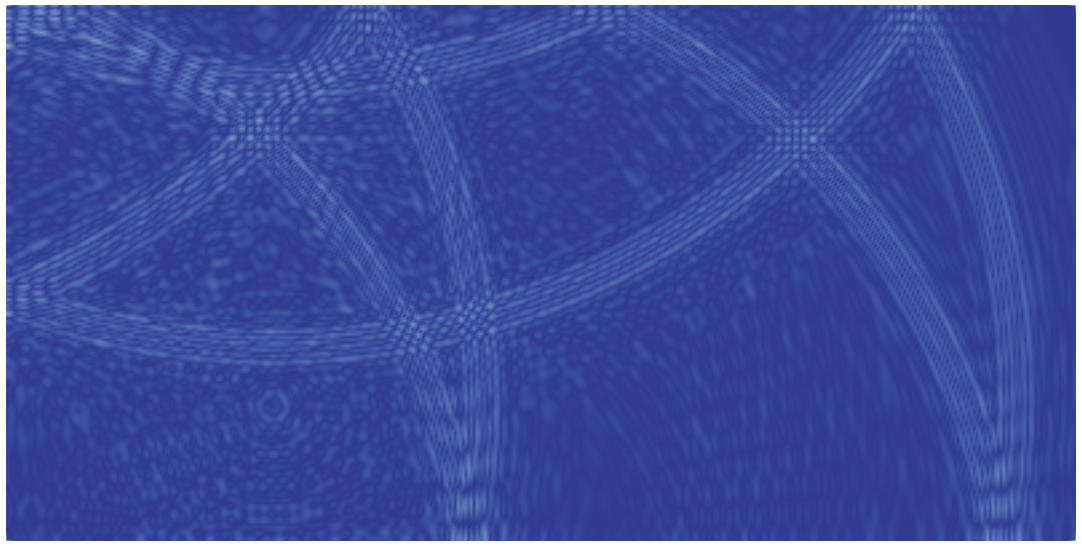}
\end{subfigure} \\[0.5em]
\begin{subfigure}[b]{0.14\textwidth}
\caption*{}
\end{subfigure}
\begin{subfigure}[b]{0.28\textwidth}
\caption{$t=3$}
\end{subfigure}
\begin{subfigure}[b]{0.28\textwidth}
\caption{$t=6$}
\end{subfigure}
\begin{subfigure}[b]{0.28\textwidth}
\centering
\caption{$t=9$}
\end{subfigure}
\caption{Absolute error in the states between the reference solution and the linear (top) and quadratic (bottom) operator inference based reduced-order models at selected time steps $t \in \{3,6,9\}$.}
\label{fig:2d_err}
\end{figure}

Because governing equation \eqref{eq:wave} contains a second-order time derivative, the reduced-order system using the quadratic manifold state approximation \eqref{eq:quadratic_approx} takes the form
\begin{equation}
\dfrac{\text{d}^2\widehat{\mathbf{s}}}{\text{d}t^2}  = \widehat{\mathbf{c}} + \widehat{\mathbf{A}} \widehat{\mathbf{s}} + \widetilde{\mathbf{H}}(\widehat{\mathbf{s}} \otimes \widehat{\mathbf{s}}); \quad \widehat{\mathbf{s}}(0) = \mathbf{V}^\top ( \mathbf{s}_0 - \mathbf{s}_\text{ref} ); \quad \dfrac{\text{d}}{\text{d}t}\widehat{\mathbf{s}}(0) = \mathbf{V}^\top \dot{\mathbf{s}}_0.
\label{eq:reduced_system_3d}
\end{equation}
The quadratic manifold formulation then leads to the following learning problem for inferring the reduced-order operators directly from the training data:
\begin{equation}
    \left(\widehat{\mathbf{c}}, \widehat{\mathbf{A}}, \widehat{\mathbf{H}} \right) = \argmin_{\substack{\widehat{\mathbf{c}} \in \mathbb{R}^r, \widehat{\mathbf{A}} \in \mathbb{R}^{r \times r}, \\ \widehat{\mathbf{H}} \in \mathbb{R}^{r \times r^2}}} \sum_{j=1}^k
 \left\|  \widehat{\mathbf{c}} + \widehat{\mathbf{A}} \widehat{\mathbf{s}}_j + \widehat{\mathbf{H}}(\widehat{\mathbf{s}}_j \otimes \widehat{\mathbf{s}}_j) - \dfrac{\text{d}^2\widehat{\mathbf{s}}_j}{\text{d}t^2} \right\|_2^2 + \lambda_1 \left( \| \widehat{\mathbf{c}} \|_2^2 + \| \widehat{\mathbf{A}} \|_F^2 \right)+ \lambda_2 \| \widehat{\mathbf{H}} \|_F^2.
 \label{eq:opinf_quadratic_second_order}
 \end{equation}
The second-order time derivatives in \eqref{eq:opinf_quadratic_second_order} are approximated from the state snapshots using a five-point central difference stencil, which is fourth-order accurate. The first two and last two time derivatives are computed using a one-sided six-point stencil to maintain the same level of accuracy.
We infer the reduced-order operators by solving linear least-squares problem \eqref{eq:opinf_quadratic_second_order}. The same procedure applies for the reduced-order models using a linear approximation subspace, the only difference being that there is no quadratic reduced matrix operator $\widehat{\mathbf{H}}$ to be sought.

The inferred reduced models are simulated with an explicit central difference time stepping scheme (which has second-order accuracy), time step size $\Delta t=10^{-2}$, and end time $T=10$. We show representative results for reduced models with a basis dimension of $r=40$, accounting for $60.1\%$ and $98.1\%$ of the snapshot energy for the linear and quadratic models, respectively. Tuning of the regularization parameters is often critical to the performance of the operator inference method. It was determined that a choice of $\lambda_1=10^{-2}$ and $\lambda_2=8.6596 \cdot 10^{-2}$ in \eqref{eq:opinf_quadratic_second_order} produces the most accurate results. For the linear reduced-order model we choose $\lambda=7.4989 \cdot 10^{-2}$.
Figure \ref{fig:2d} compares the computed solution fields for the linear and quadratic reduced-order models at selected time steps. The same snapshots of a full-order model simulation are shown for reference. The linear operator inference model does a reasonable job of capturing the wave propagation with time, as well as handling multiple reflection and interference events. However, the wave amplitudes tend to be underestimated and the model displays non-physical, oscillatory behavior throughout the computational domain. The quadratic model, on the other hand, is able to find more truthful wave amplitudes and suppresses much of the non-physical behavior away from the propagating fronts. This can be seen further in Figure~\ref{fig:trace} where we plot the operator inference predicted solution over time, monitored at the center of the computational domain.
Figure \ref{fig:2d_err} plots the absolute error in the solution at three time instances. While the error for both models is concentrated around the wave fronts, the magnitude of the absolute error is much smaller in the quadratic model.

\section{Conclusions and future work}
\label{sec:conclusions}

This paper proposes a new approach for nonlinear dimensionality reduction. The methodology learns a data-driven quadratic manifold from high-dimensional data by solving a linear regression problem. Combining the quadratic manifold state approximation with operator inference model reduction for linear problems leads to reduced-order models that are effective for problems exhibiting slow decay in the Kolmogorov \textit{n}-width. The approach can be viewed as a form of data-driven closure modeling, since the quadratic component of the manifold approximation introduces directions that lie in the orthogonal complement of the linear subspace, but does not increase the number of degrees of freedom in the low-dimensional representation. For the linear PDE examples presented in this paper, approximation in a quadratic manifold leads to a reduced-order model with quadratic structure. Using operator inference, the approach is entirely non-intrusive---that is, the manifold representation and associated reduced-order model operators can be learned directly from snapshot data.

The paper also introduces a column selection algorithm that optimally selects a subset of quadratic terms to be included in the manifold representation. This reduces the number of basis vectors required to represent the quadratic terms and reduces the dimension of the reduced-order quadratic operator to be inferred. For the examples considered in this paper, the quasi-quadratic representation incurs some error penalty and thus represents a tradeoff between reduced model performance and the computational complexity of inferring the reduced-order model.

The quadratic manifold approximation could also be used to derive reduced-order models for nonlinear problems; however, a number of computational challenges must be addressed. One set of challenges relates to the efficient representation of nonlinear terms in the reduced-order model. One way to address this is through hyperreduction, as in \cite{https://doi.org/10.48550/arxiv.2204.02462}, which introduces an additional layer of approximation beyond the state dimension reduction. Another approach is to exploit structure of the nonlinear terms through lifting to higher-order polynomial terms as in \cite{doi:10.2514/1.J057791, QIAN2020132401}, or through explicit treatment of non-polynomial nonlinear terms as in \cite{BENNER2020113433}. Determining how to combine these treatments of nonlinear dynamical systems with our quadratic manifold approximations is an important direction of future work.
We also note that the column selection method will be essential when applying the approach to nonlinear PDEs, since the reduced models for those problems will likely have higher-order polynomial structure. For example, using a quadratic manifold to reduce a PDE with quadratic operator structure will yield a reduced model with quartic structure. For such cases, it will be essential to reduce the number of reduced operator coefficients to be inferred, by using column selection to identify multi-way interactions between reduced-order state variables that can be neglected.
A second set of challenges in extending the approach to nonlinear problems relates to the algorithms to determine the basis $\overline{\mathbf{V}}$. With increased complexity of dynamics in a nonlinear system, the dimension of the reduced space is expected to grow, and it may become increasingly difficult to maintain a well conditioned least-squares problem. Exploring alternative numerical formulations, including other strategies for regularization, is another important direction of future work.

\section*{Acknowledgments}

We thank the members of the Willcox research group at the Oden Institute for their insights and many useful discussions. This work has been supported in part by the U.S.\ Department of Energy AEOLUS MMICC center under award DE-SC0019303, program manager W.\ Spotz, and by the AFOSR MURI on physics-based machine learning, award FA9550-21-1-0084, program manager F.\ Fahroo.

\appendix
\section{The SpaRSA algorithm}
\label{app:sparse}

We now discuss a basic proximal-gradient algorithm for solving the sum-of-$\ell_2$ regularized problem \eqref{eq:least-squares_reg2}. This is essentially the SpaRSA algorithm described in \cite{wright2009sparse}. The basic step of the algorithm moves from the current estimate $\overline{\mathbf{V}}$ to a new estimate $\overline{\mathbf{V}}_+$ by solving the following subproblem for some step length $\alpha>0$:
\begin{equation} \label{eq:sparsa.1}
    \overline{\mathbf{V}}_+ = \argmin_{\overline{\mathbf{V}}_+ \in \mathbb{R}^{n \times r(r+1)/2} }  \left( \frac{1}{2\alpha}  \left\| \overline{\mathbf{V}}_+^\top - \left[ \overline{\mathbf{V}}^\top - \alpha \mathbf{W} (\mathbf{W}^\top \overline{\mathbf{V}}^\top - \boldsymbol{\mathcal{E}} ) \right] \right\|_F^2 + \lambda \sum_{j=1}^{r(r+1)/2} \| \overline{\mathbf{V}}^\top_{+,j} \|_2 \right),
\end{equation}
where $\overline{\mathbf{V}}^\top_{+,j}$ denotes the $j$th row of $\overline{\mathbf{V}}^\top_+$. In fact, we can separate the objective in \eqref{eq:sparsa.1} according to rows of $\overline{\mathbf{V}}_+^\top$, and write it as
\begin{equation} \label{eq:sparsa.2}
    \overline{\mathbf{V}}_+ = \argmin_{\overline{\mathbf{V}}_+ \in \mathbb{R}^{n \times r(r+1)/2} }  \sum_{j=1}^{r(r+1)/2} \left( \dfrac{1}{2\alpha}  \left\| \overline{\mathbf{V}}_{+,j}^\top - \left[ \overline{\mathbf{V}}^\top - \alpha \mathbf{W} (\mathbf{W}^\top \overline{\mathbf{V}}^\top - \boldsymbol{\mathcal{E}}  ) \right]_{j,.} \right\|_F^2 + \lambda  \| \overline{\mathbf{V}}^\top_{+,j} \|_2 \right).
\end{equation}
Defining $\mathbf{R}_{\alpha}:= \overline{\mathbf{V}}^\top - \alpha \mathbf{W} (\mathbf{W}^\top \overline{\mathbf{V}}^\top - \boldsymbol{\mathcal{E}} ) $ and $\mathbf{R}_{\alpha,j}$ to be the $j$th row of $\mathbf{R}_{\alpha}$, we can rewrite \eqref{eq:sparsa.2} as
\begin{equation} \label{eq:sparsa.3}
     \overline{\mathbf{V}}_+ = \argmin_{\overline{\mathbf{V}}_+ \in \mathbb{R}^{n \times r(r+1)/2} }  \sum_{j=1}^{r(r+1)/2} \left( \dfrac{1}{2\alpha}  \left\| \overline{\mathbf{V}}_{+,j}^\top - \mathbf{R}_{\alpha,j} \right\|_F^2 + \lambda  \| \overline{\mathbf{V}}^\top_{+,j} \|_2 \right).
\end{equation}
The closed-form solution of  \eqref{eq:sparsa.3} is
\begin{equation}
    \label{eq:sparsa.4}
    \overline{\mathbf{V}}^\top_{+,j} = \begin{cases}
     0 & \;\; \mbox{if } \| \mathbf{R}_{\alpha,j} \|_2 \le \alpha \lambda \\
     \left(1-\dfrac{\alpha \lambda}{\| \mathbf{R}_{\alpha,j}\|_2} \right) \mathbf{R}_{\alpha,j}, & \;\; \mbox{otherwise}
    \end{cases}, \;\; j=1,2,\dotsc, r(r+1)/2.
\end{equation}

\begin{algorithm}
\caption{SpaRSA Algorithm for solving \eqref{eq:least-squares_reg2}}
\label{alg:sparsa}
\begin{algorithmic}
\STATE Given $\lambda \ge 0$, Choose $\bar{\alpha}>0$, $\overline{\mathbf{V}}_0$; \\
\STATE $\alpha \leftarrow \bar\alpha$; \\
\FOR{$k=0,1,2,\dotsc$}
\STATE Evaluate $\mathbf{R}_{\alpha}$ with $\overline{\mathbf{V}}=\overline{\mathbf{V}}_k$ and solve \eqref{eq:sparsa.4} for $\overline{\mathbf{V}}_+$; \\
\WHILE{$\mathcal{F}_{\lambda}(\overline{\mathbf{V}}_+) \ge \mathcal{F}_{\lambda}(\overline{\mathbf{V}}_k)$}
\STATE $\alpha \leftarrow \alpha/2$;
\STATE Evaluate $\mathbf{R}_{\alpha}$ with $\overline{\mathbf{V}}=\overline{\mathbf{V}}_k$ and solve \eqref{eq:sparsa.4} for $\overline{\mathbf{V}}_+$;
\ENDWHILE
\STATE Set $\overline{\mathbf{V}}_{k+1} \leftarrow \overline{\mathbf{V}}_+$ and $\alpha \leftarrow \min( (3/2) \alpha, \bar\alpha)$;
\IF{$k>5$}
\IF{$\dfrac{\mathcal{F}_{\lambda}(\overline{\mathbf{V}}_k) - \mathcal{F}_{\lambda} (\overline{\mathbf{V}}_{k-5})}{\mathcal{F}_{\lambda}(\overline{\mathbf{V}}_{k-5})} \le \epsilon_{\text{tol}}$}
\STATE Terminate;
\ENDIF
\ENDIF
\ENDFOR
\end{algorithmic}
\end{algorithm}

In most contexts, this algorithm is not applied for a single value of the regularization parameter $\lambda$ but rather a range of distinct values. The main loop in Algorithm~\ref{alg:sparsa} can be enclosed in an outer loop which iterates over these $\lambda$ values in decreasing order, with the solution $\overline{\mathbf{V}}$ for one value of $\lambda$ being used as the starting point for the next smaller value of $\lambda$. The overall procedure is detailed in \cite{wright2009sparse}. From among these different solutions, some external criterion can be used to select the most desirable or appropriate. In the current context, we might target a certain number $q$ of nonzero columns in $\overline{\mathbf{V}}$.

\bibliography{bib}
\bibliographystyle{ieeetr}

\end{document}